\documentclass{amsart}
\usepackage{amssymb}

\newtheorem{theorem}{Theorem}[section]
\newtheorem{proposition}{Proposition}[section]
\newtheorem{corollary}{Corollary}[section]
\newtheorem{lemma}{Lemma}[section]

\newtheorem*{factA}{Fact A}
\newtheorem*{factB}{Fact B}
\newtheorem*{factC}{Fact C}
\newtheorem*{KC}{Kaplansky's Conjecture}
\newtheorem*{claim}{\sl Claim}
\newtheorem*{claimone}{\sl Claim 1}
\newtheorem*{claimtwo}{\sl Claim 2}
\newtheorem*{claimthree}{\sl Claim 3}
\newtheorem*{claimfour}{\sl Claim 4}
\newtheorem*{claimfive}{\sl Claim 5}

\theoremstyle{definition}
\newtheorem*{definition}{Definition}
\newtheorem*{definitions}{Definitions}
\newtheorem*{comment}{Comment}
\newtheorem*{comments}{Comments}
\newtheorem*{notations}{Notations}
\newtheorem{remark}{Remark}[section]
\newtheorem{example}{Example}[section]
\newtheorem*{notation}{Notation}
\newtheorem*{convention}{Convention}

\newcommand{\dontwrite}[1]{}

\newcommand{\twone}{\mbox{$\text{\rm II}_1$}}
\DeclareMathOperator*{\wlim}{\text{\sf W}\text{\rm -lim}}
\begin{document}
\title[Spectral Symmetry in \twone-Factors]{Spectral Symmetry in
II${}_{\mathbf{1}}$-Factors}
\author{Sang Hyun Kim}
\address{Department of Mathematics, Kansas State university, Manhattan KS 66506, U.S.A.}
\email{shk@math.ksu.edu}
\author{Gabriel Nagy}
\address{Department of Mathematics, Kansas State university, Manhattan KS 66506, U.S.A.}
\email{nagy@math.ksu.edu}
\begin{abstract}
A self-adjoint element in a finite AW*-factor is
{\em spectrally symmetric}, if its spectral measure under the
quasitrace is invariant under the change of variables $t\longmapsto -t$.
We show that if $\mathcal{A}$ is an AW*-factor of type \twone, a self-djoint
element in $\mathcal{A}$, without full support, has quasitrace zero, if and
only if it can be written as a sum of at most three commuting spectrally symmetric
elements.
\dontwrite{
Moreover, an arbitray self-adjoint
element in $\mathcal{A}$ has quasitrace zero if and only if it can be written
as a can be written -- in $\text{\rm Mat}_2(\mathcal{A})$ -- as a
sum of at most three commuting spectrally symmetric
elements.
}
\end{abstract}
\maketitle

\section*{Introduction}

According to the Murray-von Neumann classification, finite
von Neumann factors are either of type I${}_{\text{fin}}$, or of type
\twone. For the non-expert, the easiest way to understand this classification is
by accepting the famous result of Murray and von Neumann (see \cite{MvN})
which states that {\em every finite von Neumann factor $\mathcal{M}$
 posesses a  unique state-trace $\tau_{\mathcal{M}}$}.
Upon accepting this result, the type of $\mathcal{M}$ is decided by so-called
{\em dimension range\/}:
$\mathcal{D}_{\mathcal{M}}=\big\{\tau_{\mathcal{M}}(P)\,:\,P\text{ projection in
}\mathcal{M}\big\}$ as follows. If $\mathcal{D}_{\mathcal{M}}$ is finite, then
$\mathcal{M}$ is of type I${}_{\text{fin}}$ (more explictly, in this case
$\mathcal{D}_{\mathcal{M}}=\big\{\frac kn\,:\,k=0,1,\dots,n\big\}$ for some
$n\in\mathbb{N}$, and $\mathcal{M}\simeq \text{Mat}_n(\mathbb{C})$ -- the 
algebra of $n\times n$ matrices). If $\mathcal{D}_{\mathcal{M}}$ is infinite, then
$\mathcal{M}$ is of type \twone, and in fact one has $\mathcal{D}_{\mathcal{M}}=[0,1]$.
From this point of view, the factors of type \twone\ are the ones
that are interesting, one reason being the fact that, although all factors of type
\twone\ have
the same dimension range, there are uncountably many non-isomorphic ones
(by a celebrated result of Connes).

In this paper we deal with a very simple problem. We start with a von Neumann
\twone-factor $\mathcal{M}$, a (self-adjoint)
element $A\in\mathcal{M}$, and we wish to characterize the
condition: $\tau_{\mathcal{M}}(A)=0$. The main feature of the trace $\tau_{\mathcal{M}}$
is
\begin{equation}
\tau_{\mathcal{M}}(XY-YX)=0,\,\,\,\forall\,X,Y\in\mathcal{M},
\label{tauXY-YX=0}
\end{equation}
so a sufficient condition for $\tau_{\mathcal{M}}(A)=0$
is that $A$ be expressed as a {\em sum of commutators}, i.e.
of elements of the form $[X,Y]=XY-YX$ with $X,Y\in\mathcal{M}$.
A remarkable result due to Fack and de la Harpe (\cite{FdH}) states not only that
this condition is sufficient, but if $A=A^*$ then $A$ can be written as a
sum of at most five commutators of the form $[X,X^*]$.

The aim of this paper is to characterize the condition $\tau_{\mathcal{M}}(A)=0$
in a way that is ``Hilbert space free.'' What we have in mind of
course is the purely algebraic setting due to Kaplansky (\cite{Kap}), who attempted to
formalize the theory of von Neumann algebras without any use of pre-duals.
What emerged from Kaplansky's work was the concept of AW*-algebras, which we
recall below.
\begin{definition}
A unital C*-algebra $\mathcal{A}$ is called an {\em AW*-algebra}, if for every non-empty
set $\mathcal{X}\subset\mathcal{A}$, the left anihilator
set
$\mathbf{L}(\mathcal{X})=\big\{A\in\mathcal{A}\,:\,AX=0,\,\,\,
\forall\,X\in\mathcal{X}\big\}$
is the principal right ideal generated by a projection $P\in\mathcal{A}$, that is,
 $\mathbf{L}(\mathcal{X})=\mathcal{A}P$.
\end{definition}
One can classify the finite AW*-factors into the types
I${}_{\text{fin}}$ and \twone, exactly as above, but using the
following alternative result: {\em any finite AW*-factor $\mathcal{A}$ posesses
a unique normalized quasitrace $q_{\mathcal{A}}$}. Recall that a
{\em quasitrace\/} on a C*-algebra $\mathfrak{A}$ is a map $q:\mathfrak{A}\to\mathbb{C}$ with the
following properties:
\begin{itemize}
\item[(i)] if $A,B\in\mathfrak{A}$ are self-adjoint, then
$q(A+iB)=q(A)+iq(B)$;
\item[(ii)] $q(AA^*)=q(A^*A)\geq 0$, $\forall\,A\in\mathfrak{A}$;
\item[(iii)] $q$ is linear on all abelian C*-subalgebras of $\mathfrak{A}$,
\item[(iv)] there is a map $q_2: \text{Mat}_2(\mathcal{A})\to\mathbb{C}$ with
properties (i)-(iii), such that
$$
q_2\left(\left[
\begin{array}{cc}
A& 0\\ 0 &0
\end{array}
\right]\right)=q(A),
\,\,\,\forall\,A\in\mathfrak{A}.
$$
\end{itemize}
(The condition that $q$ is {\em normalized\/} means that $q(I)=1$.)

With this terminology, the dimension range of a finite AW*-factor is the set
$\mathcal{D}_{\mathcal{A}}=\big\{q_{\mathcal{A}}(P)\,:\,P\text{ projection in
}\mathcal{A}\big\}$, and the classification into the two types
is eaxctly as above. As in the case of von Neumann factors, one can show that
the AW*-factors of type I${}_{\text{fin}}$ are again the matrix algebras
$\text{Mat}_n(\mathbb{C})$, $n\in\mathbb{N}$. The type \twone\ case however is
still mysterious. In fact, a longstanding problem in the theory of AW*-algebras
is the following:
\begin{KC}
Every AW*-factor of type \twone\ is a von Neumann factor.
\end{KC}
An equivalent formulation states that: {\em if $\mathcal{A}$ is an AW*-factor of type
\twone, then the quasitrace $q_{\mathcal{A}}$ is linear\/}
(so it is in fact a trace).
Why does one restrict Kaplansky's Conjecture to the
case of factors? On the one hand, as Dixmier had shown, there are
examples of abelian AW*-algebras which are not von Neumann algebras. Such
algebras are those that lack the existence of sufficiently many normal states.
(The spectra of such algebras are precisely those {\em stonean\/} spaces that
are not {\em hyperstonean}.) On the other hand however,
one has the following.

\begin{factA}
If $\mathcal{B}$ is an abelian AW*-subalgebra of a finite AW*-factor
$\mathcal{A}$, then $\mathcal{B}$ is a von Neumann algebra.
\end{factA}

\noindent (This is due to the fact that the restriction $\tau=q_{\mathcal{A}}\big|_{\mathcal{B}}:\mathcal{B}\to\mathbb{C}$
is normal and faithful.)
Actually (see \cite{Ped} for example), a bit more can be said, namely:
\begin{factB}
If $\mathcal{M}$ is an AW*-subalgebra of a finite AW*-factor
$\mathcal{A}$, such that the restriction
$q_{\mathcal{A}}\big|_{\mathcal{M}}\to\mathbb{C}$ is linear,
then $\mathcal{M}$ is a von Neumann algebra.
\end{factB}

A remarkable result of Haagerup (\cite{Ha}) states that {\em if $\mathcal{A}$
is an AW*-factor
of type \twone, generated (as an AW*-algebra) by an exact C*-algebra, then
$\mathcal{A}$ is a von Neumann algebra}. A simple application of this
result gives:
\begin{factC}
Any
AW*-factor of
type \twone
contains a unital
AW*-subalgebra $\mathcal{R}$ 
that
is $*$-isomorphic to
the hyperfinite
von Neumann \twone-factor.
\end{factC}
\noindent
(This can be proven exactly as in the von Neumann case, buliding up
first a copy of the diadic UHF algebra, and taking the AW*-completion.)

Suppose now $\mathcal{A}$ is an AW*-factor of type \twone, and
$A\in\mathcal{A}$ is a self-adjoint element with $q_{\mathcal{A}}(A)=0$.
In our search for a characterization of this condition, it is worth pointing
out that, in the von Neumann case, the elements of the form $XX^*-X^*X$ are the
ones that are ``certain to have trace zero,'' whereas in the AW*-factor setting,
this is not known to be the case. The natural question that arises in
connection with this
observation is: {\em which self-adjoint elements in $\mathcal{A}$ are
``certain to have quasitrace zero''?} Since the only subsitute for
\eqref{tauXY-YX=0} is
\begin{equation}
q_{\mathcal{A}}(UBU^*)=q_{\mathcal{A}}(B),\,\,\,\forall\,B\in\mathcal{A},\,
U\in\mathbf{U}(\mathcal{A}),
\label{q-uinv}
\end{equation}
where $\mathbf{U}(\mathcal{A})$ denotes the unitary group of $\mathcal{A}$,
our supply of such elements, can consist of  those
self-adjoint elements $B\in\mathcal{A}$, for which
there exists a unitary $U\in\mathbf{U}(\mathcal{A})$ with $UBU^*=-B$.
It turns out that
one can go even beyond these elements, by considering those self-adjoint $B$'s
which are {\em spectrally symmetric in $\mathcal{A}$}.
This notion will be made precise in
Section 2, but roughly speaking it means that the positive eigenvalues are
the same as the negative eigenvalues, with equal multiplicities (which are
computed using the quasitraces of the spectral projections).
Using this (still imprecise) terminology, the main result of this paper
states that
{\em a self-adjoint element $A$ with $q_{\mathcal{A}}(A)=0$ can be written,
after a suitable matrix stabilization, as the sum of
three commuting spectrally symmetric elements}.

The paper is organized as follows.
In Section 1 we discuss a certain type of
convergence for nets in AW*-factors of type \twone, which is adequate when
dealing with abelian ones.
Section 2 covers the basic properties of approximate unitary
equivalence and spectral symmetry.
Section 3 deals with a certain  integration
technique that is inspired from von Neumann's minimax trace formula (see
\cite{MvN} and \cite{FdH}).
Section 4 contains most of the technical results.
The main results are
containd in Section 5.

Parts of this paper overlap with the first author's PhD dissertation. The first author
wishes to express his gratitude to his thesis
advisor Gabriel Nagy, for essential contributions to this project.

\

\section{Weak convergence}

In this section we discuss a possible substitute for weak convergence
in the AW*-setting. We begin by adopting the following terminology.
Given an AW*-algebra $\mathcal{A}$, we call a subalgebra $\mathcal{M}\subset
\mathcal{A}$ a {\em von Neumann subalgebra}, if
\begin{itemize}
\item $\mathcal{M}$ is an AW*-algebra of $\mathcal{A}$;
\item $\mathcal{M}$ is a von Neumann algebra, i.e. $\mathcal{M}$ is a dual
Banach space.
\end{itemize}
The starting point in our discussion is the observation that {\em
AW*-subalgebras of von Neuman subalgebras are von Neumann subalgebras}.

\begin{definition}
Let $\mathcal{A}$ be an AW*-algebra. We say that a net
$(A_\lambda)_{\lambda\in\Lambda}\subset\mathcal{A}$ is {\em weakly convergent in
$\mathcal{A}$}, if
if there exists a von Neumann subalgebra $\mathcal{M}$ of $\mathcal{A}$, such that
\begin{itemize}
\item[(i)]  there exists some $\lambda_{\mathcal{M}}\in\Lambda$,
such that $A_\lambda\in\mathcal{M}$, $\forall\,\lambda\succ\lambda_{\mathcal{M}}$;
\item[(iii)] the net $(A_\lambda)_{\lambda\in\Lambda}$ is convergent
in $\mathcal{M}$ in the $\text{\sc w}^*_{\mathcal{M}}$-topology.
\end{itemize}
Observe that in this case, if $\mathcal{N}$ is any
other von Neumann subalgebra of $\mathcal{A}$ with property (i), then
it will satisfy condition
(ii) automatically. Indeed, if we choose
$\mu\in\Lambda$ such that $\mu\succ\lambda_{\mathcal{M}}$ and
$\mu\succ\lambda_{\mathcal{N}}$, then $A_\lambda\in\mathcal{M}\cap\mathcal{N}$,
$\forall\,\lambda\succ\mu$. Moreover $\mathcal{M}\cap \mathcal{N}$ is a von
Neumann subalgebra
in both $\mathcal{M}$ and $\mathcal{N}$, so
one will have the equalities
$$\text{\sc w}^*_{\mathcal{M}}\lim_{\lambda\in\Lambda}A_\lambda=
\text{\sc w}^*_{\mathcal{N}}\lim_{\lambda\in\Lambda}A_\lambda=
\text{\sc w}^*_{\mathcal{M}\cap\mathcal{N}}
\lim_{\lambda\in\Lambda}A_\lambda\in\mathcal{M}\cap\mathcal{N}.$$
In particular, the limit $\text{\sc w}^*_{\mathcal{M}}
\lim_{\lambda\in\Lambda}A_\lambda$ is independent on the particular choice of $\mathcal{M}$
-- as long as $\mathcal{M}$ has properties (i)-(ii). This element will then be
denoted by $\wlim_{\lambda\in\Lambda}A_\lambda$, and will
be referred to as the weak limit of the net $(A_\lambda)_{\lambda\in\Lambda}$.
(When there is any danger of confusion, the notation
$\text{\sf W}_{\mathcal{A}}\text{-lim}$ will be used.)
\end{definition}

\begin{remark}
Assume $\mathcal{A}$ is an AW*-algebra, and $\mathcal{B}$ is an AW*-subalgebra
of $\mathcal{A}$. For a net $(A_\lambda)_{\lambda\in\Lambda}\subset\mathcal{B}$,
the conditions:
\begin{itemize}
\item[(i)] $(A_\lambda)_{\lambda\in\Lambda}$ is weakly convergent in
$\mathcal{A}$, and
\item[(ii)] $(A_\lambda)_{\lambda\in\Lambda}$ is weakly convergent in
$\mathcal{B}$,
\end{itemize}
are equivalent, and moreover one has the equality
$\text{\sf W}_{\mathcal{A}}\text{-lim}_{\lambda\in\Lambda}A_\lambda=
\text{\sf W}_{\mathcal{B}}\text{-lim}_{\lambda\in\Lambda}A_\lambda$.

Indeed, if condition (i) is satisfied, there exists some von Neumann
subalgebra $\mathcal{M}\subset\mathcal{A}$, such that
$A_\lambda\in\mathcal{M}$, $\forall\,\lambda\succ\lambda_M$, and
some element $A\in\mathcal{M}$, such that $\text{\sc w}^*_{\mathcal{M}}
\lim_{\lambda\in\Lambda}A_\lambda=A$. In this case, we simply notice that
$\mathcal{N}=\mathcal{M}\cap\mathcal{B}$ is a von Neumann subalgebra of
$\mathcal{B}$ (hence also of $\mathcal{A}$),
so by the above discussion we must have
$\text{\sc w}^*_{\mathcal{N}}
\lim_{\lambda\in\Lambda}A_\lambda=A$.
The implication (ii)$\Rightarrow$ (i) is trivial,
since any von Neumann subalgebra of $\mathcal{B}$ is also
a von Neumann subalgebra of $\mathcal{A}$.
\end{remark}

\begin{comment}
In what follows we are going to restrict ourselves to the case when the
ambient AW*-algebra $\mathcal{A}$ is a {\em finite factor}. In this case the key
observation is the fact that (see the introduction)
if $\mathcal{M}\subset\mathcal{A}$ is an AW*-algebra with the property that
the restriction $q_{\mathcal{A}}\big|_{\mathcal{M}}:\mathcal{M}\to\mathbb{C}$
is linear, then $\mathcal{M}$ is a von Neumann algebra. In particular, {\em
all abelian AW*-subalgebras of $\mathcal{A}$ are von Neumann subalgebras}.
\end{comment}

With the above discussion in mind, the
following terminology will be useful.

\begin{definitions}
Let $\mathcal{A}$ be a $*$-algebra. For 
$\mathcal{X}\subset\mathcal{A}$,
define $\mathcal{X}^*=\{X^*\,:\,X\in\mathcal{X}\}$.
\begin{itemize}
\item[(a)] A subset $\mathcal{X}\subset\mathcal{A}$ is said to be
{\em abelian}\, if $XY=YX$, $\forall\,X,Y\in\mathcal{X}$.

\item[(b)] A subset $\mathcal{X}\subset\mathcal{A}$ is said to be
{\em $*$-abelian}\, if $\mathcal{X}\cup\mathcal{X}^*$ is abelian.

\item[(c)] A subset $\mathcal{X}\subset\mathcal{A}$ is said to be
{\em involutive}\, if $\mathcal{X}^*=\mathcal{X}$.
\end{itemize}
It is obvious that, if $\mathcal{X}$ is involutive, then
``$*$-abelian'' is equivalent to ``abelian.'' This is the case for instance
when $\mathcal{X}\subset\mathcal{A}_{sa}(=\{A\in\mathcal{A}\,:\,A=A^*\})$.
\end{definitions}

\begin{remark}
If $\mathcal{A}$ is a finite AW*-factor, and if $\mathcal{X}\subset\mathcal{A}$ is
a $*$-abelian subset, then $\mathcal{X}$ is contained in an abelian
von Neumann subalgebra $\mathcal{M}\subset\mathcal{A}$, for example
$\mathcal{M}=(\mathcal{X}\cup\mathcal{X}^*)''$ -- the bicommutant
of $\mathcal{X}\cup\mathcal{X}^*$ in $\mathcal{A}$. 
\end{remark}

\begin{remark}
The above observation is useful when dealing with
$*$-abelian nets. More explicitly, if $\mathcal{A}$ is a finite
AW*-factor, and $(A_\lambda)_{\lambda\in\Lambda}\subset\mathcal{A}$ is
a $*$-abelian net -- as a set $\mathcal{X}=\{A_\lambda\,:\,\lambda\in\Lambda\}$ -- then the condition, that
$(A_\lambda)_{\lambda\in\Lambda}$ is weakly convergent in $\mathcal{A}$, is equivalent
to the condition that $(A_\lambda)_{\lambda\in\Lambda}$ is convergent
in $\mathcal{M}=(\mathcal{X}\cup\mathcal{X}^*)''$ in the
$\text{\sc w}^*_{\mathcal{M}}$-topology.
\end{remark}

\begin{remark}
If $\mathcal{A}$ is a finite AW*-factor, 
the operation of taking
weak limits of abelian nets in $\mathcal{A}$ is ``quasi-linear,'' in the
following sense.
If $(A_\lambda)_{\lambda\in\Lambda}$ and
$(B_\lambda)_{\lambda\in\Lambda}$ are weakly convergent
{\em jointly\/} $*$-abelian nets in
$\mathcal{A}$, meaning that the set
$$\mathcal{X}=\{A_\lambda\,:\,\lambda\in\Lambda\}\cup
\{B_\lambda\,:\,\lambda\in\Lambda\})$$ is $*$-abelian,
then for any $\zeta\in\mathbb{C}$, the (abelian) net $(A_\lambda +\zeta
B_\lambda)_{\lambda\in\Lambda}$ is weakly convergent, with limit
$$\wlim_{\lambda\in\Lambda}(A_\lambda +\zeta B_\lambda)
=\big[\wlim_{\lambda\in\Lambda}A_\lambda\big]
+\zeta\big[\wlim_{\lambda\in\Lambda}A_\lambda\big].$$
\end{remark}

\begin{lemma}
Let $\mathcal{A}$ be a finite AW*-factor,
and let
$(A_\lambda)_{\lambda\in\Lambda}$ be an abelian
net in
$\mathcal{A}_{sa}$, which is
\begin{itemize}
\item bounded, i.e. $\sup_{\lambda\in\Lambda}\|A_\lambda\|<\infty$, and
\item monotone, i.e. has one of the properties $(\uparrow)$ or
$(\downarrow)$ below:
\begin{itemize}
\item[$(\uparrow)$] $\lambda_1\succ\lambda_2\Rightarrow A_{\lambda_1}
\geq A_{\lambda_2}$,
\item[$(\downarrow)$] $\lambda_1\succ\lambda_2\Rightarrow A_{\lambda_1}\leq A_{\lambda_2}$.
\end{itemize}
\end{itemize}
Then the net $(A_\lambda)_{\lambda\in\Lambda}$ is weakly convergent. Moreover,
if we take $A=\wlim_{\lambda\in\Lambda}A_\lambda$, then for any integer
$k\geq 1$, one has the following properties:
\begin{itemize}
\item[(i)] the net $(A^k_\lambda)_{\lambda\in\Lambda}$ is weakly convergent,
and $\wlim_{\lambda\in\Lambda}A_\lambda^k=A^k$;
\item[(ii)] $q_{\mathcal{A}}(A^k)=\lim_{\lambda\in\Lambda}q_{\mathcal{A}}(A_\lambda^k)$.
\end{itemize}
\end{lemma}
\begin{proof}
Consider the bicommutant
$\mathcal{M}=\{A_\lambda\,:\,\lambda\in\Lambda\}''$, which is a
von Neumann algebra.
By Remark
1.3, in order to prove the first statement, and statement (i),
it suffices to show that: the nets
$(A_\lambda^k)_{\lambda\in\Lambda}$, $k\in\mathbb{N}$ are all $\text{\sc
w}^*_{\mathcal{M}}$ convergent, and moreover,
$\text{\sc w}^*_{\mathcal{M}}\lim_{\lambda\in\Lambda}A_\lambda=A^k$,
$\forall\,k\in\mathbb{N}$, where $A=\text{\sc w}^*_{\mathcal{M}}\lim_{\lambda\in\Lambda}A_\lambda$.
The fact that the nets $(A^k_\lambda)_{\lambda\in\Lambda}$,
$k\in\mathbb{N}$ are convergent is clear, since all these nets are monotone
and bounded (the fact that $(A_\lambda)_{\lambda\in\Lambda}$ is abelian is key
for the monotonicity).
To prove the second assertion,
we
define
$X_k=\text{\sc w}^*_{\mathcal{M}}\lim_{\lambda\in\Lambda}A_\lambda^k$,
and we notice that, due to the monotonicity and boundedness of the
nets $(A^k_\lambda)_{\lambda\in\Lambda}$, we actually have:
$$X_k=\text{\sc
so-}\lim_{\lambda\in\Lambda}A^k_\lambda\text{ (in $\mathcal{M}$)},$$
where ``{\sc so}'' stands for the {\em
strong operator topology}, (coming from a realization of $\mathcal{M}$ as
a von Neumann algebra on some Hilbert space). Since
$(A_\lambda)_{\lambda\in\Lambda}$ is bounded, this gives
$$\text{\sc so-}\lim_{\lambda\in\Lambda}A_\lambda^k=A^k\text{ (in
$\mathcal{M}$)},$$
so we indeed have the equalities $X_k=A^k$.
Finally, since
$q_{\mathcal{A}}\big|_{\mathcal{M}}$ is a normal linear functional,
it follows that
$$\lim_{\lambda\in\Lambda}q_{\mathcal{A}}(A_\lambda^k)=
q_{\mathcal{A}}(A^k),\,\,\,\forall\,k\geq 1.\qedhere$$
\end{proof}

\

\section{Approximate Unitary Equivalence and Spectral Symmetry}

\begin{notations}
Let $\mathcal{A}$ be a unital C*-algebra.

\noindent{\bf A.}
Two elements $A, B\in\mathcal{A}$ are said to be {\em orthogonal},
in which case we write $A\perp B$, if:
$AB=BA=AB^*=B^*A=0$.
(Using the Fuglede-Putnam Theorem, in the case when one of the two is
normal, the above condition reduces to: $AB=BA=0$. If both $A$
and $B$ are normal, one only needs $AB=0$.)
A collection $(A_j)_{j\in J}\subset\mathcal{A}$ is said to be orthogonal, if
$A_i\perp A_j$, $\forall\,i\neq j$.

\noindent{\bf B.} We denote by
$\mathbf{U}(\mathcal{A})$
the group of unitaries in $\mathcal{A}$.
We denote by $\mathbf{P}(\mathcal{A})$ the collection of
projections in $\mathcal{A}$, that is,
$\mathbf{P}(\mathcal{A})=\big\{P\in\mathcal{A}_{sa}\,:\,P=P^2\big\}$.
\end{notations}

\begin{definition}
Let $\mathcal{A}$ be a unital C*-algebra. Two elements $A,B\in\mathcal{A}$ are
said to be {\em approximately unitarily equivalent}, if there exists a
sequence $(U_n)_{n=1}^\infty\subset\mathbf{U}(\mathcal{A})$ such that
$\lim_{n\to\infty}\|U_nAU_n^*-B\|=0$. In this case we write
$A\sim B$.
\end{definition}
The following result (perhaps well known) collects several easy properties of $\sim$.
\begin{proposition}
Let $\mathcal{A}$ be a unital C*-algebra.
\begin{itemize}
\item[(i)]
The relation $\sim$ is an equivalence relation
on $\mathcal{A}$.
\item[(ii)] If $A=\lim_{n\to\infty}A_n$ and $B=\lim_{n\to\infty}B_n$ (in norm),
and if $A_n\sim B_n$, $\forall\,n$, then $A\sim B$.
\item[(iii)] If $A,B\in\mathcal{A}$ are such that
$A\sim B$, then $A$ and $B$ have the same norm and the same spectrum.
\item[(iv)] If $A,B\in\mathcal{A}$ are such that
$A\sim B$, then $A^*\sim B^*$.
\item[(v)] If $A,B\in\mathcal{A}$ are such that $A\sim B$ and $A$ is normal,
then $B$ is also normal,
and furthermore $f(A)\sim f(B)$ for any continuous function $f:\mathbb{C}\to\mathbb{C}$.
\end{itemize}
\end{proposition}
\begin{proof}
(i). The reflexivity is trivial. The symmetry is clear because of the equality
$$\|UAU^*-B\|=\|U(A-U^*BU)U^*\|=\|A-U^*BU\|,\,\,\,\forall\,A,B\in\mathcal{A},
\,U\in\mathbf{U}(\mathcal{A}).$$
The transitivity is a consequence of the inequality:
\begin{gather*}
\|VUAU^*V^*-C\|\leq\|V(UAU^*-B)V^*\|+\|VBV^*-C\|\\
=\|UAU^*-B\|+\|VBV^*-C\|,\,\,\,\forall\,A,B,C\in\mathcal{A},\,U,V\in\mathbf{U}(\mathcal{A}).
\end{gather*}

(ii). This is pretty clear, since for every unitary
$U\in\mathbf{U}(\mathcal{A})$ one has the inequalities
$$\|UAU^*-B\|\leq
\|A-A_n\|+\|B-B_n\|+\|UA_nU^*-B_n\|,$$
so if we choose, for each $n$, a unitary $U_n\in\mathbf{U}(\mathcal{A})$, such that
$\|U_nA_nU_n^*-B_n\|<\frac 1n$, then $\lim_{n\to\infty}U_nAU_n^*=B$ (in norm).

(iii). Assume $A\sim B$. The equality $\|A\|=\|B\|$ is obvious, since
$\|UAU^*\|=\|A\|$, $\forall\,U\in\mathbf{U}(\mathcal{A})$.
To prove that $A$ and $B$ have the same spectrum,  it suffices (by symmetry)
to prove that, for every
$\lambda\in\mathbb{C}$, one has the implication:
$A-\lambda I\text{ invertible }\Rightarrow B-\lambda I\text{ invertible}$.
If we choose $(U_n)_{n=0}^\infty\subset\mathbf{U}(\mathcal{A})$, with
$\lim_{n\to\infty}\|U_nAU_n^*-B\|=0$, then it is trivial that
$$\lim_{n\to\infty}\|U_n(A-\lambda I)U_n^*-(B-\lambda I)\|=0,$$
so $B-\lambda I$
is the (norm) limit of a sequence $X_n=U_n(A-\lambda I)U_n^*$,
$n\geq 0$,
whose terms are all invertible elements. Since
$\|X_n^{-1}\|=\|(A-\lambda I)^{-1}\|$, $\forall\,n\geq 0$, we get
$$\lim_{n\to\infty}\|I-X_n^{-1}(B-\lambda I)\|=
\lim_{n\to\infty}\|I-(B-\lambda I)X_n^{-1}\|=0,$$
so for $n$ large both $X_n^{-1}(B-\lambda I)$ and $(B-\lambda I)X_n^{-1}$ are
invertible, and so is $B-\lambda I$.

(iv). This is trivial, since
$$\|UA^*U^*-B^*\|=\|(UAU^*-B)^*\|=\|UAU^*-B\|,\,\,\,
\forall\,A,B\in\mathcal{A},\,U\in\mathbf{U}(\mathcal{A}).$$

(v). Assume $A$ is normal, and $A\sim B$. If we choose a sequence of unitaries
$(U_n)_{n=0}^\infty\subset\mathbf{U}(\mathcal{A})$ with $B=\lim_{n\to\infty}
U_nAU_n^*$ (in norm), then by (iv) we also have
$B^*=\lim_{n\to\infty}
U_nA^*U_n^*$ (in norm), so we get the equalities
\begin{align*}
BB^*&=\lim_{n\to\infty}(U_nAU_n^*)(U_nA^*U_n^*)=\lim_{n\to\infty}U_nAA^*U_n^*,\\
B^*B&=\lim_{n\to\infty}(U_nA^*U_n^*)(U_nAU_n^*)=\lim_{n\to\infty}U_nA^*AU_n^*,
\end{align*}
(in norm) so we clearly have $BB^*=B^*B$.
Notice now that we also have
$$
B^kB^*{}^\ell=\lim_{n\to\infty}U_nA^kA^*{}^\ell U_n^*,\,\,\,\forall\,k,\ell\geq 0,$$
so in fact we get
$$p(B,B^*)=\lim_{n\to\infty}U_np(A,A^*)U_n^*,$$
for every polynomial $p(t,s)$ of two variables. Using the Stone-Weierstrass
Theorem, one then immediately gets $f(B)=\lim_{n\to\infty}U_nf(A)U_n^*$, for
every continuous function $f:K\to\mathbb{C}$, where $K$ denotes the
spectrum of $A$ (which is the same as the spectrum of $B$).
\end{proof}

Below we take a closer look at approximate unitary equivalence,
in the case when
ambient C*-algebra is a finite
AW*-factor.
To make matters a bit simpler,
we restrict our attention to self-adjoint elements.

\begin{notation}
Let $\mathcal{A}$ be a finite AW*-factor.
The restriction
of the quasitrace $q_{\mathcal{A}}$ to $\mathbf{P}(\mathcal{A})$ will be denoted by
$D_{\mathcal{A}}$ (or simply $D$, when there is no danger of confusion). The map
$D:\mathbf{P}(\mathcal{A})\to[0,1]$ is referred to as the
{\em dimension function\/} on $\mathcal{A}$.
\end{notation}

For future reference, we collect
the important properties of
the dimension function, in the following.

\begin{proposition}
Let $\mathcal{A}$ be a finite
 AW*-factor.
\begin{itemize}
\item[(i)] For $P,Q\in\mathbf{P}(\mathcal{A})$, the following are equivalent:
\begin{itemize}
\item[$\bullet$] $P\sim Q$;
\item[$\bullet$] there exists $U\in\mathbf{U}(\mathcal{A})$, such that
$UPU^*=Q$;
\item[$\bullet$] there exists $V\in\mathcal{A}$ with $VV^*=P$ and $V^*V=Q$;
\item[$\bullet$] $D(P)=D(Q)$.
\end{itemize}
\item[(ii)] 
If a collection
$(P_j)_{j\in J}\subset\mathbf{P}(\mathcal{A})$
is orthogonal, 
then
$$D\big(\bigvee_{j\in J}P_j\big)=\sum_{j\in J}D(P_j).$$
\end{itemize}
\end{proposition}
\begin{proof}
See \cite{Kap}.
\end{proof}

\begin{notations}
Let $\mathcal{A}$ be a finite AW*-factor.

\noindent{\bf A.} For an element
$A\in\mathcal{A}_{sa}$, we denote by $\mu_{\mathcal{A}}^A$ the spectral measure of
$A$ under the quasitrace $q_{\mathcal{A}}$. (If there is no danger of confusion,
we are going to omit the subscript $\mathcal{A}$ from the notation.)
To define rigourously $\mu^A$, we have to consider the space
$C_0(\mathbb{R})$ of all continuous complex-valued
functions on $\mathbb{R}$, which vanish at $\pm\infty$, and we use Riesz' Theorem
to define $\mu^A$ to be the unique (probability) measure on $Bor(\mathbb{R})$ --
the Borel $\sigma$-algebra -- which satisfies the equality
$$q_{\mathcal{A}}\big(f(A)\big)=\int_{\mathbb{R}}f\,d\mu^A,\,\,\,\forall\,f\in C_0(\mathbb{R}).$$
The measure $\mu^A$ will be called the {\em scalar spectral measure of $A$,
relative to $\mathcal{A}$}.

\noindent{\bf B.} Given an 
element $A\in\mathcal{A}_{sa}$, its
bicommutant $\{A\}''$ is an
abelian von Neumann algebra (by the discussion in Section 1).
For any Borel set $B\subset\mathbb{R}$ we denote by
$e_B:\mathbb{R}\to\mathbb{R}$ its indicator function, and then using Borel
functional calculus in $\{A\}''$ we can
construct a projection, denoted $e_B(A)\in\{A\}''$.
By construction, one has the equality
\begin{equation}
D\big(e_B(A)\big)=\mu^A(B),\,\,\,\forall\,B\in Bor(\mathbb{R}).
\label{dimeB}
\end{equation}
\end{notations}

With these notations, one has the following result.
\begin{theorem}
Let $\mathcal{A}$ be a finite 
AW*-factor.
For two
elements $A,B\in\mathcal{A}_{sa}$, the following are equivalent:
\begin{itemize}
\item[(i)] $A\sim B$;
\item[(ii)] $\mu^A=\mu^B$,
as measures on $Bor(\mathbb{R})$;
\item[(iii)] $D\big(e_{(-\infty,\lambda)}(A)\big)=
D\big(e_{(-\infty,\lambda)}(B)\big)$, $\forall\,\lambda\in\mathbb{R}$;
\item[(iv)] $D\big(e_{(-\infty,\lambda]}(A)\big)=
D\big(e_{(-\infty,\lambda]}(B)\big)$, $\forall\,\lambda\in\mathbb{R}$;
\item[(v)] $q_{\mathcal{A}}\big(f(A)\big)=
q_{\mathcal{A}}\big(f(B)\big)$, for every continuous
function $f:\mathbb{R}\to\mathbb{C}$;
\item[(vi)] $q_{\mathcal{A}}(A^k)=q_{\mathcal{A}}(B^k)$, $\forall\,k\in\mathbb{N}$.
\end{itemize}
\end{theorem}
\begin{proof}
(i)$\Rightarrow$ (vi). By Proposition 2.1, it suffices to consider the
case $k=1$. Assume $A\sim B$, so there exists $(U_n)_{n=0}^\infty
\subset\mathbf{U}(\mathcal{A})$ such that $B=\lim_{n\to\infty}U_nAU_n^*$
(in norm). Since the quasitrace is norm continuous (see \cite{BH}), the equality
$q_{\mathcal{A}}(A)=q_{\mathcal{A}}(B)$ follows immediately from
\eqref{q-uinv}.

(vi)$\Rightarrow$ (v). Assume (vi), 
fix a continuous function $f:\mathbb{R}\to\mathbb{C}$, and
let us prove the equality
\begin{equation}
q_{\mathcal{A}}\big(f(A)\big)=
q_{\mathcal{A}}\big(f(B)\big).
\label{qfA=B}
\end{equation}
Using Stone-Weierstrass Theorem, and the norm continuity of
$q_{\mathcal{A}}$, it suffices to prove \eqref{qfA=B} in the case
when $f$ is a polynomial function. (Indeed, if we consider the
compact set
$\Omega=\text{Spec}(A)\cup\text{Spec}(B)$, then
$f(A)$ and $f(B)$ depend only on
the restriction $f\big|_\Omega$, and if we choose a sequence
$(p_n)_{n=0}^\infty$ of
polynomials in one variable, such that $f\big|_\Omega=\lim_{n\to\infty}p_n$
in $C(\Omega)$, then $f(A)=\lim_{n\to\infty}p_n(A)$ and
$f(B)=\lim_{n\to\infty}p_n(B)$, in norm. Using the norm continuity of
$q_{\mathcal{A}}$ we have
$q_{\mathcal{A}}\big(f(A)\big)=\lim_{n\to\infty}
q_{\mathcal{A}}\big(p_n(A)\big)$
and
$q_{\mathcal{A}}\big(f(B)\big)=\lim_{n\to\infty}
q_{\mathcal{A}}\big(p_n(B)\big)$.)
%
When $f$ is a polynomial function however, the equality
\eqref{qfA=B} follows immediately from (vi)
combined with the linearity of $q_{\mathcal{A}}$ on each of the abelian
C*-subalgebras
$C^*(\{I,A\})$ and $C^*(\{I,B\})$.

(v)$\Rightarrow$ (ii).
Condition (v) implies
$$\int_{\mathbb{R}}f\,d\mu^A=
\int_{\mathbb{R}}f\,d\mu^B,\,\,\,\forall\,f\in C_0(\mathbb{R}),$$
so it will clearly force $\mu^A=\mu^B$.

(ii)$\Rightarrow$ (iii) and (ii)$\Rightarrow$ (iv) are
trivial because the conditions (iii) and (iv) read:
\begin{itemize}
\item[(iii)] $\mu^A\big((-\infty,\lambda)\big)=
\mu^B\big((-\infty,\lambda)\big)$, $\forall\,\lambda\in\mathbb{R}$;
\item[(iv)] $\mu^A\big((-\infty,\lambda]\big)=
\mu^B\big((-\infty,\lambda]\big)$, $\forall\,\lambda\in\mathbb{R}$.
\end{itemize}
The same argument shows that we also have the implications
(iii)$\Rightarrow$ (ii) and (iv)$\Rightarrow$ (ii), the reason being the
fact that if one considers the collections
\begin{align*}
\mathcal{J}_1&=\big\{(-\infty,\lambda)\,:\,\lambda\in\mathbb{R}\big\},\\
\mathcal{J}_2&=\big\{(-\infty,\lambda]\,:\,\lambda\in\mathbb{R}\big\},
\end{align*}
then for $k=1,2$ one has:
\begin{itemize}
\item[$(a)$] $Bor(\mathbb{R})=\mathbf{S}(\mathcal{J}_k)$ -- the $\sigma$-ring generated
by $\mathcal{J}_k$,
\item[$(b)$] $J,K\in\mathcal{J}_k\Rightarrow J\cap K\in\mathcal{J}_k$,
\end{itemize}
and then by standard arguments one has the implication
$$\mu^A\big|_{\mathcal{J}_k}=\mu^B\big|_{\mathcal{J}_k}\Rightarrow \mu^A=\mu^B.$$
(It is key here that both $\mu^A$ and $\mu^B$ are probability measures.)

(iii)$\Rightarrow$ (i). Assume condition (iii). Replacing $A$ with
$\delta A+\lambda I$, and $B$ with $\delta B+\lambda I$, with
$\delta,\lambda\in\mathbb{R}\smallsetminus\{0\}$ suitably chosen
(use also Proposition 2.1), we
can assume that $0\leq A,B\leq \alpha I$ for some $\alpha\in (0,1)$, that is,
\begin{equation}
\text{Spec}(A)\cup\text{Spec}(B)\subset[0,1).
\label{specin01}
\end{equation}
For every integer $n\geq 1$, consider then the spectral projections
$$P_{kn}=e_{[\frac{k-1}n,\frac{k}n)}(A)\text{ and }
Q_{kn}=e_{[\frac{k-1}n,\frac{k}n)}(B),\,\,\,k=1,\cdots,n,$$
which have dimensions
\begin{align*}
D(P_{kn})&=D\big(e_{(-\infty,\frac{k}n)}(A)\big)-
D\big(e_{(-\infty,\frac{k-1}n)}(A)\big),\\
D(Q_{kn})&=
D\big(e_{(-\infty,\frac{k}n)}(B)\big)-
D\big(e_{(-\infty,\frac{k-1}n)}(B)\big).
\end{align*}
Using the hypothesis (iii) we get
$D(P_{kn})=D(Q_{kn})$,
so there exist partial isometries $V_{kn}\in\mathcal{A}$ with
$V_{kn}V_{kn}^*=P_{kn}$ and $V_{kn}^*V_{kn}=Q_{kn}$.
By \eqref{specin01} we also have the
equalities
$$P_{1n}+P_{2n}+\dots +P_{nn}=Q_{1n}+Q_{2n}+\dots+Q_{nn}=I,\,\,\,\forall\,n\geq 2,$$
and then the element $U_n=\sum_{k=1}^nV_{kn}^*$ will be a unitary, satisfying
\begin{equation}
U_nP_{kn}U_n^*=Q_{kn},\,\,\,\forall\,n\geq k\geq 1.
\label{up=q}
\end{equation}
Using \eqref{specin01}, for every $n\geq 1$, one has the inequalities
\begin{align*}
\sum_{k=1}^n\tfrac{k-1}nP_{kn}&\leq A\leq
\sum_{k=1}^n\tfrac knP_{kn},\\
\sum_{k=1}^n\tfrac{k-1}nQ_{kn}&\leq B\leq
\sum_{k=1}^n\tfrac knQ_{kn}.
\end{align*}
In particular, the elements
$A_n=\sum_{k=1}^n\frac knP_{kn}$ and $B_n=\sum_{k=1}^n\frac knQ_{kn}$, will
satisfy
\begin{equation}
\|A_n-A\|\leq\tfrac 1n\text{ and }
\|B_n-B\|\leq\tfrac 1n,\,\,\,\forall\,n\geq 1,
\label{AnA}
\end{equation}
as well as:
$$U_nA_nU_n^*=B_n,\,\,\,\forall\,n\geq 1.$$
Using \eqref{AnA} we have
$$
\|U_nAU_n^*-B\|\leq\|U_nAU_n^*-U_nA_nU_n^*\|+\|B_n-B\|\leq\tfrac 2n,\,\,\,
\forall\,n\geq 1,$$
so $A$ and $B$ are indeed approximatively unitarily equivalent.
\end{proof}

\begin{corollary}
Let $\mathcal{A}$ be a finite
AW*-factor,
and let
$A_1,A_2,B_1,B_2\in\mathcal{A}_{sa}$ be 
elements,
with $A_1\sim A_2$, $A_1\perp B_1$, and $A_2\perp B_2$. The following are equivalent:
\begin{itemize}
\item[(i)] $B_1\sim B_2$;
\item[(ii)] $A_1+B_1\sim A_2 +B_2$.
\end{itemize}
\end{corollary}
\begin{proof}
Denote for simplicity $A_1+B_1$ by $X_1$ and $A_2+B_2$ by $X_2$. Using the orthogonality
assumptions, one has the equalities
$X_1^k=A_1^k+B_1^k$ and $X_2^k=A_2^k+B_2^k$,
which in turn imply the equalities
$$q_{\mathcal{A}}(X_1^k)=
q_{\mathcal{A}}(A_1^k)+
q_{\mathcal{A}}(B_1^k)\text{ and }
q_{\mathcal{A}}(X_2^k)=q_{\mathcal{A}}(A_2^k)+
q_{\mathcal{A}}(B_2^k),\,\,\,\forall\,k\in\mathbb{N}.$$
Since we have $q_{\mathcal{A}}(A_1^k)=q_{\mathcal{A}}(A_2^k)$, it follows that the
conditions
\begin{itemize}
\item[(i$'$)] $q_{\mathcal{A}}(B_1^k)=q_{\mathcal{A}}(B_2^k)$, $\forall\,k\in\mathbb{N}$,
\item[(ii$'$)] $q_{\mathcal{A}}(X_1^k)=q_{\mathcal{A}}(X_2^k)$, $\forall\,k\in\mathbb{N}$,
\end{itemize}
are equivalent. By Theorem 2.1 however we have the equivalences
(i)$\,\Leftrightarrow\,$(i$'$) and
(ii)$\,\Leftrightarrow\,$(ii$'$).
\end{proof}

The following result is a slight (but useful) improvement
of part (ii) from Proposition 2.1.

\begin{proposition}
Let $\mathcal{A}$ be a finite 
AW*-factor.
Assume
$(A_\lambda)_{\lambda\in\Lambda}$ and $(B_\lambda)_{\lambda\in\Lambda}$ are
nets in $\mathcal{A}_{sa}$,
indexed by the
same directed set $\Lambda$.
Assume that:
\begin{itemize}
\item each of the nets $(A_\lambda)_{\lambda\in\Lambda}$ and
$(B_\lambda)_{\lambda\in\Lambda}$ is bounded, abelian, and monotone;
\item $A_\lambda\sim B_\lambda$, $\forall\,\lambda\in\Lambda$.
\end{itemize}
If we consider the
self-adjoint elements
$A=\wlim_{\lambda\in\Lambda}A_\lambda$ and $B=\wlim_{\lambda\in\Lambda}B_\lambda$
(which exist by Lemma 1.1), then $A\sim B$.
\end{proposition}
\begin{proof}
Using Lemma 1.1, one has the equalities
$$q_{\mathcal{A}}(A^k)=\lim_{\lambda\in\Lambda}q_{\mathcal{A}}(A^k_\lambda)
\text{ and }
q_{\mathcal{A}}(B^k)=\lim_{\lambda\in\Lambda}q_{\mathcal{A}}(B^k_\lambda),
\,\,\,\forall\,k\in\mathbb{N}.$$
Using the hypothesis $A_\lambda\sim B_\lambda$, and Theorem 2.1,
we know that we have
$q_{\mathcal{A}}(A^k_\lambda)=q_{\mathcal{A}}(B^k_\lambda)$,
$\forall\,\lambda\in\Lambda$, $k\in\mathbb{N}$, so the above
equalities give
$$q_{\mathcal{A}}(A^k)=
q_{\mathcal{A}}(B^k),
\,\,\,\forall\,k\in\mathbb{N},$$
and the desired conclusion follows again from Theorem 2.1.
\end{proof}

In preparation for the next definition, we introduce the
following.

\begin{notations}
Let $\mathcal{A}$ be a C*-algebra.
For $A\in\mathcal{A}_{sa}$, we denote by
$A^+$ and $A^-$ the positive and negative parts of $A$ respectively.
Recall that $A^\pm\in\mathcal{A}$ are two
uniquely determined positive elements in $\mathcal{A}$, such that
$A=A^+-A^-$ and $A^+\perp A^-$.
(In fact $A^\pm =f^\pm(A)$, where $f^\pm:\mathbb{R}\to [0,\infty)$ are
the continuous functions defined by
$f^+(t)=\max\{t,0\}$ and $f^-(t)=\max\{-t,0\}$, $\forall\,t\in\mathbb{R}$.)
\dontwrite{
Remark that $C^*\{A\}=C^*\{A^+,A^-\}$, and
\begin{equation}
A^k=(A^+)^k+(-1)^k(A^-)^k,\,\,\,\forall\,k\in\mathbb{N}.
\label{powerA=}
\end{equation}
}
\end{notations}

We now introduce the main concept used in this paper.
\begin{definition}
Let $\mathcal{A}$ be a unital C*-algebra.
An element $A\in\mathcal{A}_{sa}$ is said to be {\em spectrally symmetric in
$\mathcal{A}$}, if its positive and negative parts $A^+$ and $A^-$ are
approximately unitarily
equivalent.
\end{definition}

The following result, along the same lines as Theorem 2.1, gives several
characterizations of spectral symmetry.
\begin{theorem}
Let $\mathcal{A}$ be a finite 
AW*-factor.
For an element $A\in\mathcal{A}$, the following are equivalent:
\begin{itemize}
\item[(i)] $A$ is spectrally symmetric;
\item[(ii)] $A\sim -A$;
\item[(iii)] the map $\phi:\mathbb{R}\ni t\longmapsto -t\in\mathbb{R}$
leaves the scalar spectral measure $\mu^A$  invariant,
that is, $\mu^A\big(\phi(B)\big)=
\mu^A(B)$, $\forall\,B\in Bor(\mathbb{R})$;
\item[(iv)] $q_{\mathcal{A}}(A^k)=0$, for every odd non-negative integer
$k$.
\item[(v)] there exist $A_1,A_2\in\mathcal{A}_{sa}$, with
 $A_1\perp A_2$, $A_1\sim A_2$, and $A_1-A_2=A$.
\end{itemize}
\end{theorem}
\begin{proof}
%
(i)$\Rightarrow$(v). This implication is trivial, by taking $A_1=A^+$ and $A_2=A^-$.

(v)$\Rightarrow$(iv). Assume $A=A_1-A_2$, with $A_1$, $A_2$ as
in (v). Since $A_1\perp A_2$, it is pretty obvious that
$$A^k=A_1^k+(-1)^kA_2^k,\,\,\,\forall\,k\in\mathbb{N}.$$
Since by Theorem 2.1 we
also have:
$$
q_{\mathcal{A}}(A_1^k)=
q_{\mathcal{A}}(A_2^k),
\,\,\,\forall\,k\in\mathbb{N},
$$
using the linearity of restriction of
the quasitrace $q_{\mathcal{A}}$ to the abelian
C*-subalgebra $C^*\{A_1,A_2\}$,
we get
$$q_{\mathcal{A}}(A^k)=q_{\mathcal{A}}(A_1^k)-
q_{\mathcal{A}}(A_2^k)=0,$$
for every odd non-negative integer $k$.

(iv)$\Rightarrow$ (ii).
Assume condition (iv), and let us prove that $A\sim -A$. Using Theorem 2.1, it
suffices to show that $q_{\mathcal{A}}\big(A^k\big)=
q_{\mathcal{A}}\big((-A)^k\big)$,
or equivalently:
$$q_{\mathcal{A}}\big(A^k\big)=(-1)^k q_{\mathcal{A}}\big(A^k\big),
\,\,\,\forall\,k\in\mathbb{N}.$$
For even $k$, this is trivial, while for odd $k$, this follows from (iv).

(ii)$\Leftrightarrow$ (iii). This equivalence is trivial, since the measure
$\nu:Bor(\mathbb{R})\to [0,1]$ defined by $\nu(B)=\mu^A\big(\phi(B)\big)$,
$\forall\,B\in Bor(\mathbb{R})$, concides with the scalar spectral measure
$\mu^{-A}$ of the self-adjoint element $-A$. Then condition (ii) is equivalent,
by Theorem 2.1, to the equality $\mu^A=\nu$, which is precisely condition (iii).

(ii)$\Rightarrow$ (i). Assume $A\sim -A$, and let us prove that
$A^+\sim A^-$. If we consider the continuous
function $f^+:\mathbb{R}\ni t\longmapsto\max\{t,0\}\in [0,\infty)$, then by
Proposition 2.1, we know that $A^+=f^+(A)\sim f^+(-A)$. 
The desired conclusion then follows from the obvious equality $f^+(-A)=A^-$.
\end{proof}

\begin{remark}
If $A,B\in\mathcal{A}_{sa}$ are spectrally symmetric, and $A\perp B$, then
$A+B$ is spectrally symmetric.
\end{remark}

\

\section{Scales of Projections and Riemann Integration}

\begin{definition}
Let $\mathcal{A}$ be an AW*-factor of type \twone. A {\em scale of projections\/}
in $\mathcal{A}$ is a system $\mathcal{E}=(E,J)$ consisting of a sub-interval
$J\subset [0,1]$ and a map $E:J\to\mathbf{P}(\mathcal{A})$, with the following properties:
\begin{itemize}
\item[(i)] $D\big(E(t)\big)=t$, $\forall\,t\in J$;
\item[(ii)] $t<s\Rightarrow E(t)\leq E(s)$.
\end{itemize}
Depending on the various features of the interval $J$, we say that
\begin{itemize}
\item the scale $\mathcal{E}$ is {\em closed}, if $J$ is a closed interval;
\item the scale $\mathcal{E}$ is {\em full}, if $J=[0,1]$.
\end{itemize}

Occasionally, we are going to abuse the notation and denote the
collection of projections $\{E(t)\,:\,t\in J\}$ also by $\mathcal{E}$. (To avoid confusion,
when we use this notation, we are going to use the phrase ``$\mathcal{E}$ as a set.'')
For instance, given a projection $P\in\mathbf{P}(\mathcal{A})$ we are going to use the notation
$P\in\mathcal{E}$ to indicate that $P=E(t)$, for some $t\in J$. Likewise,
if $\mathcal{P}\subset\mathbf{P}(\mathcal{A})$ is a collection of projections,
we use the notation $\mathcal{P}\subset\mathcal{E}$ to indicate that $P\in\mathcal{E}$,
$\forall\,P\in\mathcal{P}$.
\end{definition}

\begin{remark}
Let $\mathcal{A}$ be an AW*-factor of type \twone, and
let $J\subset [0,1]$ be a sub-interval.
If $\mathcal{E}=(E,J)$ is a scale over $J$, then as a set
$\mathcal{E}$ is totally ordered, and the map $E:J\to\mathcal{E}$ is
bijective, so
$J$ is equal to the set
$D(\mathcal{E})=\big\{D(E)\,:\,E\in\mathcal{E}\big\}$.
For this reason, the interval $J$ will be referred to as
the {\em dimension range of $\mathcal{E}$}.

Conversely, a totally ordered set of projections $\mathcal{E}$ is a
scale if and only if the dimension range
$D(\mathcal{E})$ is a sub-interval of
$[0,1]$.

If $\mathcal{E}=(E,J)$ is a scale of projections, and
if $J_0\subset J$ is a sub-interval, the restriction
$(E\big|_{J_0},J_0)$ is clearly a scale, which will be denoted by
$\mathcal{E}\big|_{J_0}$.
\end{remark}

\begin{remark}
If $\mathcal{A}$ is an AW*-factor of type \twone, and
$\mathcal{E}$ is a scale of projections in $\mathcal{A}$, with
dimension range $J$, then there
exists a unique closed scale $\overline{\mathcal{E}}$
with dimension $\overline{J}$
-- the closure of $J$ -- with $\overline{\mathcal{E}}\big|_J=\mathcal{E}$.
In fact, if $\overline{J}=[a,b]$, then the values at the endpoints,
which by an abuse of notation will be denoted by $E(a)$ and $E(b)$, are given
by $E(a)=\wlim_{t\to a^+}E(t)$ and
$E(b)=\wlim_{t\to b^-}E(t)$.
Because of this fact, for the remainder
of this article we are going to deal exclusively with closed scales.
The projections $E(a)$ and $E(b)$
will be referred to as the {\em initial\/} and {\em terminal\/}
projections of the scale $\mathcal{E}$.
\end{remark}

In preparation for subsequent constructions, we introduce the following:

\begin{definitions}
Let $\mathcal{A}$ be an AW*-factor of type \twone, and let $\mathcal{E}=(E,J)$
be a scale of projections in $\mathcal{A}$.

\noindent{\bf A.} Assuming the initial and terminal projections
of $\mathcal{E}$ are $F$ and $G$ respectivley,
we define
the {\em width of $\mathcal{E}$} to be the projection
$\mathbf{w}(\mathcal{E})=G-F$.
The dimension of the width, that is, the number
$D(\mathbf{w}(\mathcal{E}))$ -- which is equal to the length of the
dimension range -- will be referred
to as the {\em measure of $\mathcal{E}$}, and will be denoted by
$m(\mathcal{E})$.

\noindent{\bf B.}
Given a projection $P\in\mathcal{A}$ with
$$P\leq E(t),\,\,\,\forall\,t\in J,$$
one can define the scale $\mathcal{E}-P=(F,K)$, where
$K=J-D(P)=\big\{t-D(P)\,:\,t\in J\big\}$, and $F(t)=E(D(P)+t)-P$, $\forall\,t\in K$.
The scale $\mathcal{E}-P$ is called the {\em downward translation of
$\mathcal{E}$ by $P$}.

\noindent{\bf C.}
If $Q\in\mathcal{A}$ is a projection with $Q\perp E(t)$, $\forall\,t\in J$,
one can define the scale $\mathcal{E}+Q=(G,L)$, where
$L=J+D(Q)=\big\{t+D(Q)\,:\,t\in J\big\}$, and $G(t)=E(t-D(Q))+Q$, $\forall\,t\in L$.
The scale $\mathcal{E}+Q$ is called the {\em upward translation of
$\mathcal{E}$ by $Q$}.
\end{definitions}

\begin{comment}
Let $\mathcal{A}$ be an AW*-factor of type \twone. Given a closed scale of projections
$\mathcal{E}$ in $\mathcal{A}$, if we
translate it downward by its initial projection, we obtain a new scale,
denoted by $\tilde{\mathcal{E}}$, which has the same width, but
which has initial projection $0$. A scale with this
property is said to be {\em normalized}. Most of our subsequent constructions
will in effect depend only on the normalized scale.
\end{comment}

Scales (as sets) are characterized as follows:
\begin{proposition}
Let $\mathcal{A}$ be an AW*-factor of type \twone, and
let $J\subset [0,1]$ be a sub-interval.
Consider the collection:
$$\mathfrak{T}_{\mathcal{A}}(J)=\big\{\mathcal{P}\subset\mathbf{P}(\mathcal{A})\,:\,
\mathcal{P}\text{\rm\ totally ordered, and }D(P)\in J,\,\,\,\forall\,P\in\mathcal{P}\big\},$$
equipped with the inclusion order.
For an element $\mathcal{E}\in\mathfrak{T}_{\mathcal{A}}(J)$,
the following are equivalent:
\begin{itemize}
\item[(i)] $\mathcal{E}$ is a scale over $J$;
\item[(ii)] $\mathcal{E}$ is a maximal element in $\mathfrak{T}_{\mathcal{A}}(J)$.
\end{itemize}
\end{proposition}
\begin{proof}
(i)$\Rightarrow$ (ii). Assume $\mathcal{E}=\{E(t)\,:\,t\in J\}$ is
a scale over $J$, and let us show
that $\mathcal{E}$ is maximal in $\mathfrak{T}_{\mathcal{A}}(J)$. Start with some
$\mathcal{P}\in\mathfrak{T}_{\mathcal{A}}(J)$, with
$\mathcal{P}\supset\mathcal{E}$, and let us prove that this forces the
equality $\mathcal{P}=\mathcal{E}$. If we start with a projection
$P\in\mathcal{P}$, and if we put $t=D(P)$, then by the definition
we have $D\big(E(t)\big)=D(P)$. Put $X=E(t)-P$, and observe that,
since both $E(t)$ and $P$ are in $\mathcal{P}$,
which is totally ordered, it follows that either $X$ or $-X$ is a projection.
In either case, the equality $D\big(E(t)\big)=D(P)$
will force $D(\pm X)=0$, so we must have $X=0$, i.e. $P=E(t)$, so $P$ indeed
belongs to $\mathcal{E}$.

(ii)$\Rightarrow$ (i). Assume $\mathcal{E}$ is a maximal element in
$\mathfrak{T}_{\mathcal{A}}(J)$, and let us prove that $\mathcal{E}$ is a scale
over $J$. By Remark 3.1, all we have to prove is the equality
$D(\mathcal{E})=J$. By construction we already have $D(\mathcal{E})\subset J$,
so we only need to prove the other inclusion. We argue by the contradiction. Assume
there is some $s\in J$, such that
\begin{equation}
D(E)\neq s,\,\,\, \forall\,E\in\mathcal{E}.
\label{prop31contr}
\end{equation}
Consider the collections of projections
$$\mathcal{F}=\big\{E\in\mathcal{E}\,:\,D(E)<s\big\}
\text{ and }
\mathcal{G}=\big\{E\in\mathcal{E}\,:\,D(E)>s\big\},$$
and define the numbers
\begin{align*}
r&=\sup \big\{D(P)\,:\,P\in\{0\}\cup\mathcal{F}\big\},\\
t&=\inf \big\{D(P)\,:\,Q\in\{I\}\cup\mathcal{G}\big\}.
\end{align*}
(We add $0$ and $I$ simply because one of
$\mathcal{F}$ or $\mathcal{G}$ could be empty.)
\begin{claim}
There exist projections $P\in\{0\}\cup\mathcal{F}$ and
$Q\in\{I\}\cup \mathcal{G}$, such that
$D(P)=r$, $D(Q)=t$.
\end{claim}
To prove the existence of $P$, we may assume $r>0$.
In particular, $\mathcal{F}\neq\varnothing$, and
$r=\sup D(\mathcal{F})$.
Note that this, combined with the inequality $r\leq s$,
forces $r\in J$.
If we
consider
the set
$D(\mathcal{F})\subset[0,1]$, equipped with its natural order,
then it becomes a directed set.
For every $\lambda\in D(\mathcal{F})$ we choose
$F_\lambda\in\mathcal{F}$ with $D(F_\lambda)=\lambda$ (by total ordering of $\mathcal{F}$,
the projection $F_\lambda$ is unique), so that we get a monotone
net $(F_\lambda)_{\lambda\in D(\mathcal{F})}$ of projections.
Since we work with
projections, this forces $(F_\lambda)_{\lambda\in D(\mathcal{F})}$ to
be both abelian and
bounded, so using Lemma 1.1, this
net has a weak limit. If we put $P=\wlim_{\lambda\in D(\mathcal{F})}
F_\lambda$, then
it is obvious that $P$ is a projection, and moreover one has the
equality $D(P)=r$. Remark that, 
since $F\leq G$, $\forall\,F\in\mathcal{F},\,
G\in\mathcal{G}$ (by total ordering), one gets the inequalities
$$F\leq P\leq G,\,\,\,\forall\,F\in\mathcal{F},\,G\in\mathcal{G},$$
so the collection $\mathcal{E}\cup\{P\}$ is again totally ordered.
Notice however that, since $D(P)=r\in J$, by maximality this forces
$P\in\mathcal{E}$. Since $D(P)\leq s$, the condition \eqref{prop31contr} forces
$P\in\mathcal{F}$.
The existence of $Q$ is proven in the exact same way with the reverse order relation.

Having proven the above Claim, let us observe that, by the arguments employed in
the proof, we also have the inequalities
\begin{equation}
F\leq P\leq Q\leq G,\,\,\,\forall\,F\in\mathcal{F},\,\,G
\in\mathcal{G}.
\label{FPQG}
\end{equation}
Since we assume
\eqref{prop31contr}, it follows that $r<s<t$.
Choose then (use the properties
of AW*-algebras of type \twone) a projection $H\leq Q-P$ with
$D(H)=s-r$, and define the projection $R=P+H$, so that
$$D(R)=D(P)+D(H)=s.$$
Since we obviously have
$P\leq R \leq Q$, by \eqref{FPQG} we also get
$$F\leq R\leq G,\,\,\,\forall\,F\in\mathcal{F},\,\,G
\in\mathcal{G},$$
so the set $\mathcal{E}\cup\{R\}$ is totally ordered. Since $D(R)=s\in J$, the
maximality of $\mathcal{E}$ forces $R\in\mathcal{E}$, thus contradicting \eqref{prop31contr}.
\end{proof}

\begin{corollary}
Let $\mathcal{A}$ be an AW*-factor of type \twone, let $J\subset [0,1]$ be a
sub-interval, and let $\mathcal{P}\in\mathfrak{T}_{\mathcal{A}}(J)$. There exists
at least one scale $\mathcal{E}$ over $J$, with $\mathcal{E}\supset\mathcal{P}$.
\end{corollary}
\begin{proof}
Immediate from Zorn's Lemma, and the above characterization.
\end{proof}

\begin{remark}
Here is an interesting special case of Corollary 3.1.
Given $\mathcal{A}$ an AW*-factor of type \twone, and a self-adjoint element
$A\in\mathcal{A}$, let us consider the collection
$$\mathfrak{S}(A)=\big\{e_{(-\infty,\alpha)}(A)\,:\,\alpha\in\mathbb{R}\big\}
\cup\big\{e_{(-\infty,\beta]}(A)\,:\,\beta\in\mathbb{R}\big\}.$$
It is obvious that $\mathfrak{S}(A)$ is totally ordered. More precisely, one
has the inequalities
$$e_{(-\infty,\beta)}(A)\leq
e_{(-\infty,\beta]}(A)\leq
e_{(-\infty,\alpha)}(A),\,\,\,\forall\,\alpha>\beta.$$
Notice that $\mathfrak{S}(A)\ni 0,I$, so
by Corollary 3.1, there exists at least one full scale
$\mathcal{E}\supset\mathfrak{S}(A)$.
Such a (full) scale will be referred to as a
{\em spectral scale for $A$} (in $\mathcal{A}$).
\end{remark}

In preparation for the next construction, we introduce the following

\begin{notations}
Assume $\mathcal{A}$ is an AW*-factor of type \twone, and
$\mathcal{E}$ is a closed scale of projections
in $\mathcal{A}$, with dimension range $D(\mathcal{E})=[a,b]$.

\dontwrite{
\noindent{\bf A.} The projection $E(b)-E(a)$ will be called
the {\em support of $\mathcal{E}$}, and will be denoted
by $\mathbf{s}(\mathcal{E})$.

\noindent{\bf B.}
}
Given a partition
$$\Delta=[a=t_0<t_1<\dots<t_n=b]$$ of the
interval $[a,b]$, and a bounded function
$f:[a,b]\to\mathbb{R}$, we define the lower and upper Darboux sums
\begin{align*}
L_{\mathcal{E}}(f,\Delta)&=\sum_{k=1}^n\big[
\inf_{s\in [t_{k-1},t_k]}f(s)\big]\cdot
\big[E(t_k)-E(t_{k-1})\big],\\
U_{\mathcal{E}}(f,\Delta)&=\sum_{k=1}^n\big[
\sup_{s\in [t_{k-1},t_k]}f(s)\big]\cdot
\big[E(t_k)-E(t_{k-1})\big].
\end{align*}
Note that, for any partion $\Delta$, using the linearity of $q_{\mathcal{A}}$ on the
bicommutant $\mathcal{E}''$, one has the equalities
\begin{align}
q_{\mathcal{A}}\big(L_{\mathcal{E}}(f,\Delta)&=
\sum_{k=1}^n\big[
\inf_{s\in [t_{k-1},t_k]}f(s)\big]\cdot
[t_k-t_{k-1}]=L(f,\Delta),\label{qL=}
\\
q_{\mathcal{A}}\big(U_{\mathcal{E}}(f,\Delta)&=
\sum_{k=1}^n\big[
\sup_{s\in [t_{k-1},t_k]}f(s)\big]\cdot
[t_k-t_{k-1}]=U(f,\Delta),\label{qU=},
\end{align}
where $L(f,\Delta)$ and $U(f,\Delta)$ denote the usual (scalar)
lower and upper Darboux sums of $f$.

Observe also that, if we consider the set
$\mathfrak{P}[a,b]$ of all partitions of $[a,b]$, ordered with respect to the
inclusion, then $\mathfrak{P}[a,b]$ becomes a directed set, and moreover
\begin{itemize}
\item $\big(L_{\mathcal{E}}(f,\Delta)\big)_{\Delta\in\mathfrak{P}[a,b]}$ is an
abelian increasing net,
\item $\big(U_{\mathcal{E}}(f,\Delta)\big)_{\Delta\in\mathfrak{P}[a,b]}$ is an
abelian decreasing net.
\end{itemize}
Since we also have the inequalities
$$\big[\inf_{s\in[a,b]}f(s)\big]\cdot\mathbf{w}(\mathcal{E})
\leq
L_{\mathcal{E}}(f,\Delta)\leq
U_{\mathcal{E}}(f,\Delta)\leq
\big[\sup_{s\in[a,b]}f(s)\big]\cdot\mathbf{w}(\mathcal{E})
,\,\,\,\forall\,\Delta\in\mathfrak{P}[a,b],$$
by Lemma 1.1 these nets are weakly convergent.
\end{notations}

\begin{proposition}
With the notations above, if $f:[a,b]\to\mathbb{R}$ is Riemann integrable,
then one has the equalities
$$\wlim_{\Delta\in\mathfrak{P}[a,b]}L_{\mathcal{E}}(f,\Delta)=
\wlim_{\Delta\in\mathfrak{P}[a,b]}U_{\mathcal{E}}(f,\Delta).$$
Moreover, if we denote this common limit by $A$, then one has the equality
$$q_{\mathcal{A}}(A)=\int_a^bf(t)\,dt.$$
\end{proposition}
\begin{proof}
Put $L=\wlim_{\Delta\in\mathfrak{P}[a,b]}L_{\mathcal{E}}(f,\Delta)$ and
$U=\wlim_{\Delta\in\mathfrak{P}[a,b]}U_{\mathcal{E}}(f,\Delta)$.
If we consider the abelian von Neumann algebra
$\mathcal{M}=\mathcal{E}''$, we have the equalities
$$
L=\text{\sc w}^*_{\mathcal{M}}\lim_{\Delta\in\mathfrak{P}[a,b]}
L_{\mathcal{E}}(f,\Delta)
\text{ and }
U=\text{\sc w}^*_{\mathcal{M}}
\lim_{\Delta\in\mathfrak{P}[a,b]}U_{\mathcal{E}}(f,\Delta)
\text{ (in $\mathcal{M}$)}.$$
Using monotonicity of the nets
$\big(L_{\mathcal{E}}(f,\Delta)\big)_{\Delta\in\mathfrak{P}[a,b]}$  and
$\big(U_{\mathcal{E}}(f,\Delta)\big)_{\Delta\in\mathfrak{P}[a,b]}$, we also get
the inequalities
$$
L_{\mathcal{E}}(f,\Delta)\leq L\leq U\leq
U_{\mathcal{E}}(f,\Delta),\,\,\,\forall\,\Delta\in\mathfrak{P}[a,b].
$$
Using the order properties of the quasitrace (which is linear on $\mathcal{M}$),
we have
$$
q_{\mathcal{A}}\big(L_{\mathcal{E}}(f,\Delta)\big)\leq
q_{\mathcal{A}}(L)\leq
q_{\mathcal{A}}(U)\leq
q_{\mathcal{A}}\big(U_{\mathcal{E}}(f,\Delta)\big),
$$
which using \eqref{qL=} and \eqref{qU=} reads:
$$
L(f,\Delta)\leq
q_{\mathcal{A}}(L)\leq
q_{\mathcal{A}}(U)\leq
U(f,\Delta),
\,\,\,\forall\,\Delta,\in\mathfrak{P}[a,b].$$
Taking limit this gives
$$q_{\mathcal{A}}(L)=q_{\mathcal{A}}(U)=\int_a^bf(t)\,dt.$$
In particular
(use the linearity of the quasitrace on $\mathcal{M}$), this gives
$q_{\mathcal{A}}(U-L)=0$, and then the
inequality $U-L\geq 0$, combined with the
faitfulness of the quasitrace, will force $U=L$.
\end{proof}
\begin{notation}
Given a closed scale $\mathcal{E}$ as above -- with dimension range $[a,b]$ --
and a Riemann integrable function $f:[a,b]\to\mathbb{R}$, the element
$A\in\mathcal{E}''$, defined in the above result, will be
denoted by $\int_a^bf(t)\,dE(t)$ (or simply $\int_a^bf\,d\mathcal{E}$,
when there is no danger of confusion).
If we denote by $\mathfrak{R}[a,b]$ the algebra of real-valued Riemann integrable
functions, the correspondence
\begin{equation}
\mathfrak{R}[a,b]\ni f\longmapsto \int_a^bf\,d\mathcal{E}
\label{R-calc}
\end{equation}
will be referred to as the {\em Riemann integral calculus associated with
the scale $\mathcal{E}$}.
\end{notation}

\begin{remark}
Given a scale $\mathcal{E}$ with dimension range
$[a,b]$, and $f\in\mathfrak{R}[a,b]$, the element $A=\int_a^bf\,d\mathcal{E}$
will satisfy the inequalities
$$\big[\inf_{s\in[a,b]}f(s)\big]\cdot\mathbf{w}(\mathcal{E}) 
\leq A\leq
\big[\inf_{s\in[a,b]}f(s)\big]\cdot\mathbf{w}(\mathcal{E}).$$
(This follows from the corresponding inequalities for lower and upper Darboux
sums, after taking weak limit in $\mathcal{E}''$.)
This will then give the inequality
$$\mathbf{s}(A)\leq\mathbf{w}(\mathcal{E}),$$
where $\mathbf{s}(A)$ denotes the {\em support of $A$}.
(Recall that, given an AW*-algebra $\mathcal{A}$ and an element $A\in\mathcal{A}_{sa}$,
one defines $\mathbf{s}(A)=I-P$, where
$P\in\mathbf{P}(\mathcal{A})$ is the projection defined by the
condition $\mathbf{L}\big(\{A\})=\mathcal{A}P$.
Equivalently, using Borel functional calculus, $\mathbf{s}(A)=
e_{\mathbb{R}\smallsetminus\{0\}}(A)$.)
\end{remark}

The following technical result deals with sequential approximation.

\begin{lemma}
Let $\mathcal{A}$ be an AW*-factor of type \twone, let $\mathcal{E}$ be a scale
of projections in $\mathcal{A}$ with dimension range $[a,b]$, and let $f\in\mathfrak{R}[a,b]$.
Assume $(\Delta_n)_{n=1}^\infty$ is a sequence of partitions of $[a,b]$, with
$\Delta_1\subset \Delta_2\subset\dots$.
\begin{itemize}
\item[(i)] If $\int_a^bf(t)\,dt=\lim_{n\to\infty}L(f,\Delta_n)$, then
$\int_a^bf\,d\mathcal{E}=\wlim_{n\to\infty}L_{\mathcal{E}}(f,\Delta_n)$.
\item[(ii)] If $\int_a^bf(t)\,dt=\lim_{n\to\infty}U(f,\Delta_n)$, then
$\int_a^bf\,d\mathcal{E}=\wlim_{n\to\infty}U_{\mathcal{E}}(f,\Delta_n)$.
\end{itemize}
\end{lemma}
\begin{proof} It suffices to prove property (i). (To prove (ii) we simply use (i)
with $f$ replaced by $-f$.)

To prove (i), denote $\int_a^bf\,d\mathcal{E}$ simply by $B$, and define the
sequence
$(A_n)_{n=1}^\infty\subset\mathcal{E}''$ by
$A_n=L_{\mathcal{E}}(f,\Delta_n)$, $n\in\mathbb{N}$.
On the one hand,
since the sequence $(\Delta_n)_{n=1}^\infty$ is increasing, the sequence
$(A_n)_{n=1}^\infty$  is increasing.
On the other hand,
it is clear that
we have the inequalities
\begin{equation}
A_n\leq B,\,\,\, \forall\,n\in\mathbb{N}.
\label{AnleqB}
\end{equation}
Using Lemma 1.1, the
limit $A=\wlim_{n\to\infty}A_n$ exists, and it will have quasitrace
\begin{equation}
q_{\mathcal{A}}(A)=\lim_{n\to\infty}q_{\mathcal{A}}(A_n)=\lim_{n\to\infty}
L(f,\Delta_n)=\int_a^bf(t)\,dt=q_{\mathcal{A}}(B).
\label{qA=qB}
\end{equation}
Finally, working in $\mathcal{E}''$, \eqref{AnleqB} yields $A\leq B$.
By the linearity and faithfulness of $q_{\mathcal{A}}$ on $\mathcal{E}''$,
the equality \eqref{qA=qB} will force $A=B$.
\end{proof}

\begin{remark}
The construction of the element
$A=\int_a^bf\,d\mathcal{E}$ is compatible with translations.
To be more precise, if one defines
the translation maps $\Lambda_s:\mathfrak{R}[a,b]\to\mathfrak{R}[a+s,b+s]$ by
$$(\Lambda_sf)(t)=f(t-s),\,\,\,\forall\,t\in [a-s,b-s],\,f\in\mathfrak{R}[a,b],$$
then one has the following properties:
\begin{itemize}
\item[(i)] If $\mathcal{E}=(E,[a,b])$ is a scale, and $P$ is a projection
with $P\leq E(a)$, with dimension $D(P)=\delta$, then
$$
\int_a^bf\,d\mathcal{E}=\int_{a-\delta}^{b-\delta}\Lambda_{-\delta}f\,d(\mathcal{E}-P),
\,\,\,\forall\,f\in\mathfrak{R}[a,b].$$
\item[(ii)] If $\mathcal{E}=(E,[a,b])$ is a scale, and $Q$ is a projection
with $Q\perp E(b)$, with dimension $D(Q)=\delta$, then
$$
\int_a^bf\,d\mathcal{E}=\int_{a+\delta}^{b+\delta}\Lambda_{\delta}f\,
d(\mathcal{E}+Q),
\,\,\,\forall\,f\in\mathfrak{R}[a,b].$$
\end{itemize}
\dontwrite{
does not depend on the
initial projection $E(a)$ of the scale. In fact, if one defines
the left translation maps $\Lambda_s:\mathfrak{R}[a,b]\to\mathfrak{R}[a-s,b-s]$ by
$$(\Lambda_sf)(t)=f(t+s),\,\,\,\forall\,t\in [a-s,b-s],\,f\in\mathfrak{R}[a,b],$$
then one has the equality
\begin{equation}
\int_a^bf\,d\mathcal{E}=\int_0^{b-a}(\Lambda_af)\,d\tilde{\mathcal{E}}.
\label{R-int-trans}
\end{equation}
Using \eqref{R-int-trans} we may occasionally
assume that the scale we work with is normalized.
}\end{remark}

\begin{notation}
Given a scale $\mathcal{E}$ with dimension range $[a,b]$,
a function $f\in\mathfrak{R}[a,b]$, and a sub-interval $[a_1,b_1]\subset [a,b]$,
we are going to denote by $\int_{a_1}^{b_1}f\,d\mathcal{E}$ the element
$\int_a^be_{[a_1,b_1]}f\,d\mathcal{E}$. An equivalent description can be given
in terms of restriction:
$$\int_{a_1}^{b_1}\big(f\big|_{[a_1,b_1]}\big)\,d\big(\mathcal{E}\big|_{[a_1,b_1]}\big).$$
\end{notation}

The following result summarizes several easy properties
of this calculus.
\begin{proposition}
Let $\mathcal{A}$ be an AW*-factor of type \twone, and
let $\mathcal{E}$ be a closed scale of projections, with dimension range
$D(\mathcal{E})=[a,b]$.
\begin{itemize}
\item[(i)] The map \eqref{R-calc} is a real algebra homomorphism.
\item[(ii)] One has the inequality
$$\big\|\int_a^bf\,d\mathcal{E}\big\|\leq\|f\|_{\sup},\,\,\,
\forall\,f\in\mathfrak{R}[a,b],$$
where $\|\,.\,\|_{\sup}$ stands for the supremum norm.
\item[(iii)] Given a continuous function $\phi:\mathbb{R}\to\mathbb{R}$, one has
the equalities
$$\phi\big(\int_a^bf\,d\mathcal{E}\big)=
\int_a^b(\phi\circ f)\,d\mathcal{E},\,\,\,
\forall\,f\in\mathfrak{R}[a,b],$$
where the left hand side is obtained by continuos functional calculus applied to
the self-adjoint element $\int_a^bf\,d\mathcal{E}$.
\end{itemize}
\end{proposition}
\begin{proof}
(i). To prove additivity, we start with two Riemann integrable functions
$f_1,f_2:[a,b]\to\mathbb{R}$,
and we prove the equality
\begin{equation}
\int_a^b(f_1+f_2)\,d\mathcal{E}=\int_a^bf_1\,d\mathcal{E}+
\int_a^bf_2\,d\mathcal{E}.
\label{R-add}
\end{equation}
If we work in the von Neumann algebra $\mathcal{M}=
\langle\langle\mathcal{E}\rangle\rangle$,
the for every partition $\Delta\in\mathfrak{P}[a,b]$,
one obviously has the inequalities
$$L_{\mathcal{E}}(f_1,\Delta)+L_{\mathcal{E}}(f_2,\Delta)\leq
L_{\mathcal{E}}(f_1+f_2,\Delta)\leq
U_{\mathcal{E}}(f_1+f_2,\Delta)\leq
U_{\mathcal{E}}(f_1,\Delta)+U_{\mathcal{E}}(f_2,\Delta),
$$
so taking $\text{\sc w}^*_{\mathcal{M}}$-limit will give \eqref{R-add}.

The homogeneity property
$$
\int_a^b(\alpha f)\,d\mathcal{E}=\alpha\int_a^bf\,d\mathcal{E},\,\,\,
\forall\,\alpha\in\mathbb{R},\,f\in\mathfrak{R}[a,b]$$
is proven in the exact same way.

In order to prove that the correspondence
\eqref{R-calc} is multiplicative, it suffices to prove that it has the property:
\begin{equation}
\int_a^bf^k\,d\mathcal{E}=\big[\int_a^bf\,d\mathcal{E}\big]^k,\,\,\,
\forall\,f\in\mathfrak{R}[a,b],\,k\in\mathbb{N}.
\label{R-mult}
\end{equation}
Using the obvious equality
$$\int _a^b1\,d\mathcal{E}=\mathbf{w}(\mathcal{E}),$$
and the linearity, it may assume in \eqref{R-mult} that $f\geq 0$.
If we fix such an $f$, as well as
$k\in\mathbb{N}$, and we define the net
$(A_\Delta)_{\Delta\in\mathfrak{P}[a,b]}$ by
$$A_\Delta=L_{\mathcal{E}}(f,\Delta)^k,$$
then, on the one hand, by Lemma 1.1 we know that
$$\wlim_{\Delta\in\mathfrak{P}[a,b]}A_\Delta=\big[\int_a^bf\,d\mathcal{E}\big]^k.$$
On the other hand, using the fact that $f$ is non-negative, it is quite clear
that
$$A_\Delta=L_{\mathcal{E}}(f^k,\Delta),\,\,\,\forall\,\Delta\in\mathfrak{P}[a,b],$$
so we get
$$\wlim_{\Delta\in\mathfrak{P}[a,b]}A_\Delta=\int_a^bf^k\,d\mathcal{E}.$$

(ii). This inequality is trivial.

(iii). Fix $f\in\mathfrak{R}[a,b]$, as well as a continuous function
$\phi:\mathbb{R}\to\mathbb{R}$, and denote $\int_a^bf\,\mathcal{E}$ simply by $A$.
Using the Stone-Weierstrass Theorem, we know that
for every $\varepsilon>0$, there exists a polynomial function
$\phi_\varepsilon:\mathbb{R}\to\mathbb{R}$, such that
$$\big|\phi(s)-\phi_\varepsilon(s)\big|\leq\varepsilon,\,\,\,
\forall\,s\in \big[-\|f\|_{\sup},\|f\|_{\sup}\big].$$
On the one hand, by (ii) we have
$\text{Spec}(A)\subset\big[-\|f\|_{\sup},\|f\|_{\sup}\big]$, so
using the properties of
functional calculus we know that
\begin{equation}
\big\|\phi(A)-\phi_\varepsilon(A)\big\|\leq\varepsilon.
\label{R-phi-eps1}
\end{equation}
On the other hand, using (i) we know that
\begin{equation}
\int_a^b(\phi_\varepsilon\circ f)\,d\mathcal{E}=
\phi_\varepsilon(A).
\label{R-phi-eps2}
\end{equation}
Finally, since we obviously have
$$
\big|(\phi\circ f)(t)-(\phi_\varepsilon\circ f)(t)\big|\leq\varepsilon,\,\,\,
\forall\,t\in [a,b],$$
we also have the inequality
$$
\big\|\int_a^b(\phi\circ f)\,d\mathcal{E}-\int_a^b(\phi_\varepsilon\circ f)
\,d\mathcal{E}\big\|\leq\varepsilon,$$
so using \eqref{R-phi-eps1} and \eqref{R-phi-eps2} we get
$$\big\|\phi(A)-\int_a^b(\phi\circ f)\,d\mathcal{E}\big\|\leq 2\varepsilon.$$
Since this inequality is true for every $\varepsilon>0$, it forces
$$\phi(A)=\int_a^b(\phi\circ f)\,d\mathcal{E}.\qedhere$$
\end{proof}

\begin{corollary}
With the notations above,
if $\mathcal{E}$ is a scale with dimension range $[a,b]$, and if
$f,g\in\mathfrak{R}[a,b]$ are such that
$$f=g,\,\,\,\text{\rm (Lebesgue) a.e.}$$
then $\int_a^bf\,d\mathcal{E}=\int_a^bg\,d\mathcal{E}$.
\end{corollary}
\begin{proof}
If one considers
the commuting elements $X=\int_a^bf\,d\mathcal{E}$ and
$Y=\int_a^bg\,d\mathcal{E}$,
then the positive element $(X-Y)^2$ will have
quasitrace
$$q_{\mathcal{A}}\big((X-Y)^2\big)=
\int_a^b[f(t)-g(t)]^2\,dt=0,$$
which obviousy forces $X=Y$.
\end{proof}

\begin{corollary}
Let $\mathcal{A}$ be an AW*-factor of type \twone, let
$\mathcal{F}$ and $\mathcal{G}$ be closed scales of projections in
$\mathcal{A}$ with dimension ranges $[a,b]$
and $[c,d]$ respectively, and let $f:[a,b]\to\mathbb{R}$ and
$g:[c,d]\to\mathbb{R}$
be two Riemann integrable functions, such that
\begin{equation}
\int_a^bf(t)^k\,dt=\int_c^dg(s)^k\,ds,\,\,\,\forall\,k\in\mathbb{N}.
\label{cor-R-int}
\end{equation}
Then the elements $A=\int_a^bf\,d\mathcal{F}$ and $B=\int_c^d g\,d\mathcal{G}$
are approximately unitary equivalent.
\end{corollary}
\begin{proof}
By the properties of the Riemann calculus, for every $k\in\mathbb{N}$, we have
$$A^k=\int_a^bf^k\,d\mathcal{F}\text{ and }
B^k=\int_a^bg^k\,d\mathcal{G},$$
so using \eqref{cor-R-int}
we get
$$q_{\mathcal{A}}(A^k)=
q_{\mathcal{A}}(B^k),\,\,\,\forall\,k\in\mathbb{N}.$$
From Theorem 2.1 it follows that $A\sim B$.
\end{proof}

The next result
should be regarded as
a ``Change of Variable'' rule. In preparation for its formulation,
we introduce the following terminology.
\begin{notations}
Let $\mathcal{A}$ be an AW*-factor of type \twone, and let
$P\in\mathbf{P}(\mathcal{A})$ be a non-zero projection.
We define the compression
$$P\mathcal{A}P=\big\{PAP\,:\,A\in\mathcal{A}\big\}.$$
Of course, $P\mathcal{A}P$ is an AW*-subalgebra of $\mathcal{A}$,
with unit $P$, but it is also a {\em factor}, so in fact
$P\mathcal{A}P$ is itself an AW*-factor of type \twone.
Its quasitrace is then given by
\begin{equation}
q_{P\mathcal{A}P}(X)=\frac{q_{\mathcal{A}}(X)}{D_{\mathcal{A}}(P)},\,\,\,\forall\,X\in
P\mathcal{A}P,
\label{qPAP}
\end{equation}
\end{notations}

\begin{proposition}
Let $\mathcal{A}$ be an AW*-factor of type \twone, and let
$P\in\mathbf{P}(\mathcal{A})$ be a non-zero projection.
Put $\lambda=D(P)$.
\begin{itemize}
\item[(i)] If $\mathcal{E}=(E,[a,b])$ is a scale of projections in
$P\mathcal{A}P$, then the map
$E^P:[\lambda a,\lambda b]\to\mathbf{P}(\mathcal{A})$, given by
$$E^P(t)=E(t/\lambda),\,\,\,\forall\,t\in[\lambda a,\lambda b],$$
defines a scale $\mathcal{E}^P=\big(E^P ,[\lambda a,\lambda b]\big)$ in
$\mathcal{A}$. Moreover, for any $f\in\mathfrak{R}[\lambda a,\lambda b]$, one
has the equality
\begin{equation}
\int_{\lambda a}^{\lambda b}f(t)\,dE^P(t)\,\text{\rm (in $\mathcal{A}$)}\,=
\int_a^bf(\lambda t)\,dE(t)\,\text{\rm (in $P\mathcal{A}P$)}.
\label{CVPAP}\end{equation}
\item[(ii)] Conversely, if $\mathcal{F}=(F,[\alpha,\beta])$ is a scale in $\mathcal{A}$,
with $F(\beta)\leq P$, then the map
$E:[\alpha/\lambda,\beta/\lambda]\to\mathbf{P}(P\mathcal{A}P)$, given by
$$E(t)=F(\lambda t),\,\,\,\forall\,t\in [\alpha/\lambda,\beta/\lambda],$$
defines a scale $\mathcal{E}=\big(E,[\alpha/\lambda,\beta/\lambda]\big)$ in
$P\mathcal{A}P$, with $\mathcal{E}^P=\mathcal{F}$.
\end{itemize}
\end{proposition}
\begin{proof}
(i). The fact that $\mathcal{E}^P$ is a scale in $\mathcal{A}$ is quite clear,
since by \eqref{qPAP} we have
$$D_{\mathcal{A}}\big(E^P(t)\big)=
D_{\mathcal{A}}\big(E(t/\lambda)\big)=
D_{\mathcal{A}}(P)\cdot D_{P\mathcal{A}P}(E(t/\lambda)\big)=\lambda\cdot (t/\lambda)=t,
\,\,\,\forall\,t\in[\lambda a,\lambda b].$$

To prove the statement about Riemann integrals, we use the following notations:
\begin{itemize}
\item given a partition $\Delta\in\mathfrak{P}[a,b]$, say
$\Delta=[a=t_0<\dots<t_n=b]$, we define the partition
$\lambda\Delta\in\mathfrak{P}[\lambda a,\lambda b]$ as
$\lambda\Delta=[\lambda a=\lambda t_0<\dots<
\lambda t_n=\lambda b]$;
\item for a function $f:[\lambda a,\lambda b]\to\mathbb{R}$ we define
the function $f^\lambda:[a,b]\to\mathbb{R}$ by
$f^\lambda(t)=f(\lambda t)$, $\forall\,t\in [a,b]$.
\end{itemize}
With these notations, one has the following easy facts:
\begin{itemize}
\item[{\sc (a)}] the correspondence $\mathfrak{P}[a,b]\ni
\Delta\longmapsto \lambda\Delta\in\mathfrak{P}[\lambda a,\lambda b]$ is an order
preserving bijection;
\item[{\sc (b)}] for a bounded function $f:[\lambda a,\lambda b]\to\mathbb{R}$,
the function $f^\lambda:[a,b]\to\mathbb{R}$ is bounded, and, for any
sub-interval $[c,d]\subset [a,b]$, one has the equalities
$$\inf_{s\in [\lambda c,\lambda d]}f(s)=\inf_{t\in[c, d]}f^\lambda(t)
\text{\ \ and }
\sup_{s\in [\lambda c,\lambda d]}f(s)=\sup_{t\in[c,d]}f^\lambda(t).$$
\end{itemize}
Using the above two facts we see that $f\in\mathfrak{R}[\lambda a,\lambda b]
\Leftrightarrow  f\in\mathfrak{R}[a,b]$, and moreover, one has the equalities
$$L_{\mathcal{E}^P}(f,\lambda\Delta)=L_{\mathcal{E}}(f^\lambda,\Delta)
\text{ and }
U_{\mathcal{E}^P}(f,\lambda\Delta)=
U_{\mathcal{E}}(f^\lambda,\Delta),\,\,\,\forall\,\Delta\in\mathfrak{P}[a,b].$$
Taking weak limits in $\mathcal{A}$ and $P\mathcal{A}P$  respectively, then
yields the equality
$$\int_{\lambda a}^{\lambda b}f\,d\mathcal{E}^P
\,\text{(in $\mathcal{A}$)}\,=\int_a^bf^\lambda\,d\mathcal{E}
\,\text{(in $P\mathcal{A}P$)},
$$
which is precisely \eqref{CVPAP}.

(ii). This statement is trivial.
\end{proof}

The Riemann integral calculus developed above will be used in  connection with the
following key result.
\begin{theorem}
Let $\mathcal{A}$ be an AW*-factor of type \twone and let $A\in\mathcal{A}_{sa}$.
\begin{itemize}
\item[(i)] The function $\omega_A:[0,1]\to\mathbb{R}$, defined by
\begin{equation}
\omega_A(t)=\left\{
\begin{array}{cl}
\min\,\text{\rm Spec}(A)&\text{ if }t=0\\
\inf\big\{\alpha\in\mathbb{R}\,:\,D\big(e_{(-\infty,\alpha]}(A)\big)\geq
 t\big\}&\text{ if }t\in(0,1]
 \end{array}\right.
 \label{def-omega}\end{equation}
in non-decreasing, hence Riemann integrable. Moreover, one has the equality
$\omega_A(1)=\max\,\text{Spec}(A)$.
\item[(ii)] For any integer $k\geq 1$, one has the equality
$$q_{\mathcal{A}}(A^k)=\int_0^1\omega_A(t)^k\,dt.$$
\item[(iii)] For any spectral scale $\mathcal{E}$ for $A$, one has the equality
\begin{equation}
A=\int_0^1\omega_A(t)\,d\mathcal{E}(t).\label{A=int-omega}\end{equation}
\end{itemize}
\end{theorem}
\begin{proof}
The fact that $\omega_A\big|_{(0,1]}$ is non-decreasing is trivial. The equality
$\omega_A(1)=\max\,\text{Spec}(A)$ is pretty obvious, since the inequality
$D\big(e_{(-\infty,\alpha]}(A)\big)\geq 1$ is equivalent
to the equality $e_{(-\infty,\alpha]}=I$, which in turn is equivalent to $A\leq
\alpha I$.
To finish the proof of (i), we fix some $t\in (0,1]$, and
we must show that $\omega_A(0)\leq\omega_A(t)$.
We argue by contradiction, assuming $\omega_A(0)>\omega_A(t)$, so there
exists $\alpha\in\mathbb{R}$ with $\omega_A(0)>\alpha$, and
$D\big(e_{(-\infty,\alpha]}(A)\big)\geq t(>0)$. This is however impossible,
since the inequality $\alpha<\omega_A(0)$ forces $e_{(-\infty,\alpha]}(A)\leq
e_{(-\infty,\omega_A(0))}(A)$, and by construction
$e_{(-\infty,\omega_A(0))}(A)=0$.

Using Proposition 3.3, it is clear that 
property (ii) follows from property (iii).

To prove property (iii), we start off by fixing a spectral
scale (see Remark 3.3) $\mathcal{E}=\big(E(t)\big)_{t\in[0,1]}$
for $A$.
\begin{claimone}
For every $t\in [0,1]$ one has the inequalities
\begin{equation}
e_{(-\infty,\omega_A(t))}(A)\leq E(t)\leq e_{(-\infty,\omega_A(t)]}(A).
\label{thm31claim1}
\end{equation}
\end{claimone}
Since both projections
$e_{(-\infty,\omega_A(t))}(A)$ and
$e_{(-\infty,\omega_A(t)]}(A)$ belong to $\mathfrak{S}(A)\subset\mathcal{E}$, by
total ordering, all we have to prove are the corresponding
inequalities for the dimensions, i.e.
$$
D\big(e_{(-\infty,\omega_A(t))}(A)\big)\leq t\leq D\big(e_{(-\infty,\omega_A(t)]}
(A)\big),
$$
or equivalently, using the scalar spectral measure, 
\begin{equation}
\mu^A\big((-\infty,\omega_A(t))\big)\leq t\leq \mu^A\big((-\infty,\omega_A(t)]\big).
\label{minmaxcl1}
\end{equation}
Since $\omega_A(0)=\min\,\text{Spec}(A)$, we have
$\mu^A\big((-\infty,\omega_A(0))\big)=0$, so \eqref{minmaxcl1} is trivial for
$t=0$. Assume now $t\in (0,1]$.
To prove the inequalities \eqref{minmaxcl1}
we consider the non-decreasing functions
$f,g:\mathbb{R}\to\mathbb{R}$ defined by
$$f(\alpha)=\mu^A\big((-\infty,\alpha)\big)\text{ and }
g(\alpha)=\mu^A\big((-\infty,\alpha]\big),\,\,\,
\forall\,\alpha\in\mathbb{R},$$
and we consider the set $\Omega_t=\big\{\alpha\in \mathbb{R}\,:\,g(\alpha)\geq t\big\}$,
so that $\omega_A(t)=\inf\Omega_t$. Since $\mu^A$ is a measure on
$Bor(\mathbb{R})$,
we know that $g$ is continuous from the right, i.e.
$$g(\beta)=\lim_{\alpha\to\beta^+}g(\alpha),\,\,\,
\forall\,\beta\in\mathbb{R}.
$$
In particular, we have
$g\big(\omega_A(t)\big)=\lim_{\alpha\to\omega_A(t)^+}g(\alpha)\geq t$, which gives the second
inequality in \eqref{minmaxcl1}.
Since we also have
$$
f(\beta)=\lim_{\alpha\to\beta^-}g(\alpha),\,\,\,
\forall\,\beta\in\mathbb{R}.
$$
and $g(\alpha)<t$, $\forall\,\alpha<\omega_A(t)$,
we immediately get $f\big(\omega_A(t)\big)\leq t$, which is the
first inequality in \eqref{minmaxcl1}.

\begin{claimtwo}
$A\in \mathcal{E}'$, i.e. $AE(t)=E(t)A$,
$\forall\,t\in [0,1]$.
\end{claimtwo}
Fix $t\in [0,1]$, and notice that, since
$e_{(-\infty,\omega_A(t))}(A)$ commutes with $A$, it suffices to show
that $F(t)=E(t)-e_{(-\infty,\omega_A(t))}(A)$ commutes with $A$.
By Claim 1 it follows that $F(t)$ is a projection, and moreover,
$$F(t)\leq
e_{(-\infty,\omega_A(t)]}(A)-
e_{(-\infty,\omega_A(t))}(A)=
e_{\{\omega_A(t)\}}(A).$$
This obviously forces $F(t)A=\omega_A(t)F(t)=AF(t)$, and we are done.

\begin{claimthree}
For any partition $\Delta\in\mathfrak{P}[0,1]$, the lower and upper Darboux sums
of $\omega_A$ satisfy the inequalities
\begin{equation}
L_{\mathcal{E}}(\omega_A,\Delta)\leq A\leq
U_{\mathcal{E}}(\omega_A,\Delta).
\label{LU-vs-A}\end{equation}
\end{claimthree}
On the one hand, if $\Delta=(0=t_0<t_1<\dots<t_n=1)$, due to the
monotonicity of $\omega_A$, one has the equalities
\begin{align}
L_{\mathcal{E}}(\omega_A,\Delta)&=\sum_{k=1}^n\omega_A(t_{k-1})\big[E(t_k)-E(t_{k-1})\big],
\label{L-omega=}\\
U_{\mathcal{E}}(\omega_A,\Delta)&=\sum_{k=1}^n\omega_A(t_k)\big[E(t_k)-E(t_{k-1})\big],
\label{U-omega=}.
\end{align}
On the other hand, using Claim 1, we have the inequalities
$$e_{(-\infty,\omega_A(t_{k-1}))}(A)\leq
E(t_{k-1})\leq E(t_k)\leq e_{(-\infty,\omega_A(t_k)]}(A),$$
which gives the inequalities
\begin{equation}
E(t_k)-E(t_{k-1})\leq e_{[\omega_A(t_{k-1}),\omega_A(t_k)]}(A),\,\,\,
\forall\,k=1,\dots,n.
\label{E-under-e}
\end{equation}
Of course, the spectral projections $e_{[\omega_A(t_{k-1}),\omega_A(t_k)]}(A)$
satisfy the inequalities
$$
\omega_A(t_{k-1})e_{[\omega_A(t_{k-1}),\omega_A(t_k)]}(A)
\leq Ae_{[\omega_A(t_{k-1}),\omega_A(t_k)]}(A)\leq
\omega_A(t_k)e_{[\omega_A(t_{k-1}),\omega_A(t_k)]}(A),$$
so multiplying this inequality by $E(t_k)-E(t_{k-1})$, which by Claim 2 commutes
with all three sides, and using \eqref{E-under-e} we get
$$
\omega_A(t_{k-1})\big[E(t_k)-E(t_{k-1})\big]\leq
A\big[E(t_k)-E(t_{k-1})\big]\leq
\omega_A(t_k)\big[E(t_k)-E(t_{k-1})\big],$$
for all $k=1,\dots,n$.
Summing up, using the obvious equality
$$\sum_{k=1}^n\big[E(t_k)-E(t_{k-1})\big]=E(1)-E(0)=I,$$
as well as \eqref{L-omega=} and
\eqref{U-omega=}, the desired inequalities \eqref{LU-vs-A} immediately follow.

After all these preparations, we proceed with the proof of
\eqref{A=int-omega}.
First of all, we notice that by Claim 2 we know that
$\{A\}\cup\mathcal{E}$ is involutive and abelian, the
AW*-subalgebra $\mathcal{M}=\big(\{A\}\cup\mathcal{E}\big)''$
is an abelian von Neumann algebra.
Secondly, if we consider the element
$B=\int_0^1\omega_A(t)\,d\mathcal{E}(t)$, then
$A$ and $B$ belong to $\mathcal{M}$. Moreover, since one has the equalities
$$B=\text{\sc w}^*_{\mathcal{M}}\text{-}\lim_{\Delta\in\mathfrak{P}[0,1]}
L_{\mathcal{E}}(\omega_A,\Delta)=
\text{\sc w}^*_{\mathcal{M}}\text{-}\lim_{\Delta\in\mathfrak{P}[0,1]}
U_{\mathcal{E}}(\omega_A,\Delta),$$
by Claim 3 we must have both inequalities $B\leq A$ and $A\leq B$, so we
indeed have the equality $A=B$.
\end{proof}

The result below -- essentially a converse of Remark 3.4 --
is useful when estimating the dimension of the
support.

\begin{proposition}
Given  an AW*-factor $\mathcal{A}$ of type \twone, and a non-zero
element $A\in\mathcal{A}_{sa}$, there exist
\begin{itemize}
\item[(i)] a scale $\mathcal{E}=\big(E,[0,\delta])$ with
$E(\delta)=\mathbf{s}(A)$, and
\item[(ii)] a non-decreasing function $f:[0,\delta]\to
\big[\min\,\text{\rm Spec}(A),\max\,\text{\rm Spec}(A)\big]$, such that
$A=\int_0^\delta f\,d\mathcal{E}$.
\end{itemize}
\end{proposition}
(Note that (ii) in fact forces $D\big(\mathbf{s}(A)\big)=\delta$.)
\begin{proof}
Denote $\mathbf{s}(A)$ simply by $P$, and let $D(P)=\delta$. Since
$A\in P\mathcal{A}P$, one can write
$$A=\int_0^1g\,d\mathcal{F}\,\text{(in $P\mathcal{A}P$)},$$
where $\mathcal{F}=(F,[0,1])$ is a spectral scale for $A$ in $P\mathcal{A}P$
(so $F(1)=P$), and
$g\in\mathfrak{R}[0,1]$ is some non-decreasing function,
namely $\omega_A$, but computed in $P\mathcal{A}P$. Of course, since
$\text{Spec}_{P\mathcal{A}P}(A)\subset\text{Spec}(A)$, one has:
\begin{equation}
\min\,\text{Spec}(A)\leq g(t)\leq\max\,\text{Spec}(A),\,\,\,\forall\,t\in [0,1].
\label{ginSpec}
\end{equation}
Using Proposition 3.4, if we consider $\mathcal{E}=\mathcal{F}^P$, namely
$\mathcal{E}=(E,[0,\delta])$, with
$$E(t)=F(t/\delta),\,\,\,\forall\,t\in[0,\delta],$$
and if we define the function $f\in\mathfrak{R}[0,\delta]$ by
$$f(t)=g(t/\delta),\,\,\,\forall\,t\in[0,\delta],$$
then by Proposition 3.4 we get:
$$
\int_0^\delta f(t)\,dE(t)\,\text{(in $\mathcal{A}$)}\,=
\int_0^1f(\delta t)\,dF(t)\,\text{(in $P\mathcal{A}P$)}\,=
\int_0^1g\,d\mathcal{F}\,\text{(in $P\mathcal{A}P$)}\,=A.
$$
Using \eqref{ginSpec}, we also have
the inclusion $\text{Range}\,f\subset\big[\min\,\text{Spec}(A),
\max\,\text{Spec}(A)\big]$. Finally, the equality $E(\delta)=F(1)=P$
is trivial. 
\end{proof}

We conclude this section with two applications of the Riemann calculus.
One application (Proposition 3.6 below) deals with
``copying'' elements. The second one (Example 3.1) shows hwo to build
self-adjoint elements with prescribed scalar spectral measure. 

Before we discuss the next result, let us introduce the following terminology.

\begin{definition}
Let $\mathcal{A}$ be a finite AW*-factor.
A subalgebra $\mathcal{B}\subset\mathcal{A}$ is said to be
an {\em AW*-subfactor of $\mathcal{A}$}, if:
\begin{itemize}
\item[$\bullet$] $\mathcal{B}$ is an AW*-subalgebra of $\mathcal{A}$, which contains
$I$ -- the unit of $\mathcal{A}$;
\item[$\bullet$] $\mathcal{B}$ is a factor.
\end{itemize}
\end{definition}
\begin{comments} Let
$\mathcal{A}$ be a finite AW*-factor.
\begin{itemize}
\item[A.]
If $\mathcal{B}$ is an AW*-subfactor of $\mathcal{A}$, then $\mathcal{B}$ is
obviously a finite AW*-factor, and its canonical quasitrace (by uniqueness) is
given by $q_{\mathcal{B}}=q_{\mathcal{A}}\big|_{\mathcal{B}}$ -- the restriction
of $q_{\mathcal{A}}$ to $\mathcal{B}$. In particular, in the case when
$\mathcal{B}$ is of type \twone (this forces $\mathcal{A}$ to be of type
\twone as well), for a collection $\mathcal{E}\subset\mathcal{B}$, the conditions
\begin{itemize}
\item[(a)] $\mathcal{E}$ is a scale of projections in $\mathcal{A}$, and
\item[(b)] $\mathcal{E}$ is a scale of projections in $\mathcal{B}$,
\end{itemize}
are equivalent. Moreover, if $\mathcal{E}$ has dimension range
$[a,b]$, then
$$
\int_a^bf\,d\mathcal{E}\,\text{(in $\mathcal{B}$)}\,
=
\int_a^bf\,d\mathcal{E}\,\text{(in $\mathcal{A}$)},\,\,\,\forall\,
f\in\mathfrak{R}[a,b].
$$
\item[B.] It turns out (see \cite{Na} for example; this is not necessary here)
that if $\mathcal{B}$ is an arbitrary C*-subalgebra of $\mathcal{A}$, with
$\mathcal{B}\ni I$, and such that $\mathcal{B}$ is an AW*-factor (in itself), then
$\mathcal{B}$ is automatically an AW*-subfactor of $\mathcal{A}$.
\end{itemize}
\end{comments}

\begin{proposition}
Let $\mathcal{A}$ be an AW*-factor of type \twone, let
$P\in\mathbf{P}(\mathcal{A})$ be a non-zero projection, let
$A\in\mathcal{A}_{sa}$,
and let $\mathcal{B}$ be an AW*-subfactor of $P\mathcal{A}P$, of type \twone.
For any projection $Q\in\mathbf{P}(\mathcal{B})$ with $D_{\mathcal{A}}(Q)\geq
D_{\mathcal{A}}\big(\mathbf{s}(A)\big)$, there exists $B\in (Q\mathcal{B}Q)_{sa}$ with $B\sim A$.
\end{proposition}
\begin{proof}
By Proposition 3.6 there exists a scale $\mathcal{F}=(F,[0,\delta])$ in
$\mathcal{A}$, with $F(\delta)=\mathbf{s}(A)$, and
$f\in\mathfrak{R}[0,\delta]$, such that $A=\int_0^\delta f\,d\mathcal{F}$.

Denote $D_{\mathcal{A}}(P)$ simply by $\lambda$, so that
$$D_{\mathcal{B}}(Q)=
D_{P\mathcal{A}P}(Q)=\frac{D_{\mathcal{A}}(Q)}{\lambda}\geq
\frac{D_{\mathcal{A}}\big(\mathbf{s}(A)\big)}{\lambda}=\frac{\delta}{\lambda}.$$
Choose a projection $Q_0\in\mathbf{P}(\mathcal{B})$ with $Q_0\leq Q$ and
$D_{\mathcal{B}}(Q_0)=\delta/\lambda$, and let
$\mathcal{E}=(E,[0,\delta/\lambda])$ be a scale in $\mathcal{B}$ with
$E(\delta/\lambda)=Q_0$. Let us consider the element
$$B=\int_0^{\delta/\lambda} f(\lambda t)\,dE(t)\text{ (in $\mathcal{B}$)}.$$
Since $\mathcal{E}$ is also a scale in $P\mathcal{A}P$, we also have the
equality
$$B=\int_0^{\delta/\lambda} f(\lambda t)\,dE(t)\text{ (in $P\mathcal{A}P$)}.$$
Let $\mathcal{E}^P=(E^P,[0,\delta])$ be the scale in $\mathcal{A}$,
constructed in Proposition 3.4. According to Proposition 3.4, we have the
equality
$$B=\int_0^\delta f(t)\,dE^P(t)\text{ (in $\mathcal{A}$)},$$
and then by Corollary 3.3 (applied to $f=g$ and to the scales $\mathcal{F}$ and $\mathcal{E}^P$)
it follows immediately that $A\sim B$.
\end{proof}

\begin{comments}
If $\mathcal{A}$ is an AW*-factor of type \twone, the maps
$\omega_A$, associated with elements $A\in\mathcal{A}_{sa}$,
have several additional properties
listed below. (These features are not needed here; see \cite{Na2}
for details.)
\begin{itemize}
\item[A.] For any $A\in\mathcal{A}_{sa}$,
the map $\omega_A:[0,1]\to\mathbb{R}$ is continuous from the left,
continuous at $0$, and satisfies:
$$\text{Spec}(A)=
\overline{\text{Range}\,\omega_A}.$$
\item[B.] If $A\in\mathcal{A}_{sa}$ is positive, then
\begin{itemize}
\item[(i)] $\omega_A(t)=
\inf\big\{\|PAP\|\,:\,P\in\mathbf{P}(\mathcal{A}),\,\,
D(P)\geq t\big\}$, $\forall\,t\in (0,1]$;
\item[(ii)] if $\mathcal{E}=(E,[0,1])$ is a spectral scale for $A$, then
$$\omega_A(t)=\|E(t)A\|,\,\,\,\forall\,t\in (0,1].$$
\end{itemize}
\item[C.]
Given a full scale $\mathcal{E}$, and a
non-decreasing function $f:[0,1]\to\mathbb{R}$, which is continuous from the
left, and continuous at $0$,
 the element $A=\int_0^1f\,d\mathcal{E}\in\mathcal{A}_{sa}$ satisfies the identity
$\omega_A=f$. Moreover, $\mathcal{E}$ is a spectral scale for $A$.
\end{itemize}
In the case of von Neumann \twone-factors,
the map $t\longmapsto\omega_A(t)$ is related to the
singular numbers discussed in \cite{MvN}, in connection with the
min-max trace formula (which is precisely
property (ii) in Theorem 3.1 above, with $k=1$).
Using the language from \cite{MvN}, if $\mathcal{A}$ is a
von Neuman \twone-factor, and $A\in\mathcal{A}$, $A\geq 0$,
then for every $t\in [0,1]$,
one has the equality
$\omega_A(1-t)=s_t(A)$, where
$s_t(A)$ is the ``$t^{\text{th}}$ singular number of $A$.''
\end{comments}

\begin{example}
Let $\mathcal{A}$ be an AW*-factor ot type \twone, and let
$\mathcal{E}=(E,[0,1])$ be a full scale in $\mathcal{A}$. We can define the element
$M=M_{\mathcal{E}}=\int_0^1t\,dE(t)\in\mathcal{A}$. By Proposition 3.3 we know that
the scalar spectral measure $\mu^M$ of $M$ is given by
\begin{equation}
\int_{\mathbb{R}}\phi\,d\mu^M=\int_0^1\phi(t)\,dt,\,\,\,\forall\,\phi\in C_0(\mathbb{R}).
\label{med}
\end{equation}
An element $M\in\mathcal{A}_{sa}$ with property \eqref{med} is called
a {\em mediator in $\mathcal{A}$}. The specific element $M_{\mathcal{E}}$ is
referred to as the {\em mediator of $\mathcal{E}$}. It is obvious that
$\text{Spec}(M)=[0,1]$. Given a Borel measurable function
$f:[0,1]\to\mathbb{R}$, which is Riemann integrable, it is not hard to show
(see \cite{Na2}) that one has the equality
$\int_0^1f\,d\mathcal{E}=f(M_{\mathcal{E}})$.
\end{example}

We conclude with a discussion on probabilistic independence, that is necessary
in the following section.
\begin{definition}
Let $\mathcal{A}$ be an AW*-factor of type \twone, and let $\mathcal{B}$ be an
AW*-subalgebra of $\mathcal{A}$. We say that $\mathcal{B}$ is {\em thick}, if
there exist a mediator $M\in\mathcal{A}$, such that
\begin{itemize}
\item $BM=MB$, $\forall\,B\in\mathcal{B}$;
\item $q_{\mathcal{A}}(BM^k)=q_{\mathcal{A}}(B)\cdot q_{\mathcal{A}}(M^k)$,
$\forall\,B\in\mathcal{B},\,k\in\mathbb{N}$.
\end{itemize}
In this case, $M$ will be referred to as a {\em $\mathcal{B}$-mediator\/} (in
$\mathcal{A}$).
\end{definition}

Obviously the center $\mathbb{C}(=\{\lambda I\,:\,\lambda\in\mathbb{C}\})$ is thick,
and every mediator is a $\mathbb{C}$-mediator. 
Because of possible (type) limitations on the commutant, not all AW*-subalgebras
are thick. The terminology below is meant to provide a method of testing for
thickness.

\begin{definition}
Two AW*-subfactors $\mathcal{B}$ and $\mathcal{M}$ of $\mathcal{A}$ are
said to be {\em independent in probability}, if:
\begin{itemize}
\item[$\bullet$] $\mathcal{B}$ and $\mathcal{M}$ commute, i.e.
$BM=MB$, $\forall\,B\in\mathcal{B},\,M\in\mathcal{M}$;
\item[$\bullet$] $q_{\mathcal{A}}(BM)=q_{\mathcal{B}}(B)\cdot
q_{\mathcal{M}}(M)$, $\forall\,B\in\mathcal{B}_{sa},\,M\in\mathcal{M}_{sa}$.
\end{itemize}
\end{definition}

With this terminology, one has the following observation.
\begin{remark}
If $\mathcal{A}$ is an AW*-factor of type \twone, and 
$\mathcal{B}$ and $\mathcal{M}$ are two AW*-subfactors that are independent
in probability, with $\mathcal{M}$ of type \twone, then
$\mathcal{B}$ is thick in $\mathcal{A}$. In fact, every mediator $M$ in $\mathcal{M}$
(such elements exist by Example 3.1) is a $\mathcal{B}$-mediator.
\end{remark}

\begin{example}
Every AW*-factor $\mathcal{A}$, of type \twone, contains a thick
subfactor of type \twone. One way to construct such subfactors is the following.
We start with an AW*-subfactor $\mathcal{R}\subset\mathcal{A}$ that is
isomorphic to the hyperfinite von Neumann \twone-factor (see Fact C in the introduction).
Since $\mathcal{R}\simeq \mathcal{R}\otimes\mathcal{R}$ -- spatial tensor
product of von Neumann algebras -- it follows that $\mathcal{R}$ contains two
subfactors, namely $\mathcal{R}\otimes I$ and $I\otimes\mathcal{R}$, which are
obviously independent in probability.
Regarding these as AW*-subfactors of $\mathcal{A}$ finishes the
construction.
\end{example}

\

\section{Foldings}

This section consists of several technical results, necessary in Section 5.
At some point, a certain hypothesis (global for this section) will be set.

\begin{definitions}
Let $\mathcal{A}$ be an AW*-factor of type
\twone, and let $k\geq 1$ be an integer.
A double sequence $\Phi=(A_1,\dots,A_k;B_1,\dots,B_k)\subset\mathcal{A}_{sa}$ is
called a {\em $k$-folding in $\mathcal{A}$}, if:
\begin{itemize}
\item[(i)] $A_i\perp B_j$, $\forall\,i,j\in\{1,\dots,k\}$;
\item[(ii)] $\{A_1,\dots,A_k\}$ and $\{B_1,\dots,B_k\}$ are abelian.
\item[(iii)] $A_i\sim B_i$, $\forall\,i\in\{1,\dots,k\}$.
\end{itemize}
Note that, using (i), it turns out that $\{A_1,B_1,\dots, A_k,B_k\}$ is abelian.
Using (i) and (iii), it follows that the elements $S_i=A_i-B_i$, $i=1,\dots,k$,
are spectrally symmetric and they all commute.
We define the {\em support of $\Phi$} to be the projection
$$\mathbf{s}(\Phi)=\big[\bigvee_{i=1}^k\mathbf{s}(A_i)\big]\vee
\big[\bigvee_{i=1}^k\mathbf{s}(B_i)\big]$$
We define $\|\Phi\|=\max\big\{\|A_1\|,\|B_1\|,\dots,\|A_n\|,\|B_n\|\big\}$.

When we want to identify the
elements $X=A_1+\dots+A_k$ and $Y=B_1+\dots +B_k$ (which are orthogonal), we are
going to use the phrase: {\em $\Phi$ is a folding of $X$ as $Y$}. 
\end{definitions}

\begin{comment}
For $X, Y\in\mathcal{A}_{sa}$, the existence of a
folding of $X$ as $Y$, obviously implies the condition
$q_{\mathcal{A}}(X)=
q_{\mathcal{A}}(Y)$.
The main goal of this paper is essentially to prove the converse of this
statement.
\end{comment}
\begin{remark}
Suppose $\Phi_n=(A_{n1},\dots,A_{nk};B_{n1},\dots,B_{nk})$, $n=1,\dots,N$ are
$k$-foldings in $\mathcal{A}$, which are orthogonal, in the sense that
$\mathbf{s}(\Phi_m)\perp\mathbf{s}(\Phi_n)$, $\forall\,m\neq n$. Then the double sequence
$(A_1,\dots,A_k;B_1,\dots,B_k)$, defined by
$$A_j=\sum_{n=1}^NA_{nj}\text{ and }B_j=\sum_{n=1}^NB_{nj},
\,\,\,\forall\,j\in\{1,\dots,k\},$$
is a $k$-folding, which will be denoted by $\Phi_1+\dots+\Phi_n$.
This follows essentially from Corollary 2.1. It is also pretty clear that
$\|\Phi_1+\dots+\Phi_n\|=\max\big\{\|\Phi_1\|,\dots,\|\Phi_n\|\big\}$.
\end{remark}

In what follows, we are going to isolate a special type of
$2$-foldings, that consist of projections.

\begin{definitions}
A {\em superprojection in $\mathcal{A}$} is a system
$\pi=(P_1,P_2;P_3,P_4)$ of projections in $\mathcal{A}$, with the following properties:
\begin{itemize}
\item $P_i\perp P_j$, $\forall\,i\neq j$;
\item $P_1\sim P_3$ and $P_2\sim P_4$.
\end{itemize}
It is obvious that $\pi$ is a $2$-folding.

Given another superprojection $\pi'=(P_1',P_2';P_3',P_4')$, we write
$\pi'\leq \pi$, if $P_i'\leq P_i$, $i=1,2,3,4$.
\end{definitions}

\dontwrite{
\begin{convention}
For the remainder of this section we are going to work under the 
following assumptions:
{\em We fix $\mathcal{A}$ to be an AW*-factor of type \twone. We fix a thick
AW*-subfactor $\mathcal{B}$ of type \twone (which exists by Example 3.2). We fix
a $\mathcal{B}$-mediator $M$ in $\mathcal{A}$.}

For future reference we will say that $\mathcal{A}$, $\mathcal{B}$,
$\mathcal{M}$, and $M$ satisfy the {\sc Technical Condition}.
\end{convention}
}

\begin{lemma}
Let $\mathcal{A}$ be an AW*-factor of type \twone, let $\pi=(P_1,P_2;P_3,P_4)$
be a superprojection in $\mathcal{A}$, and let $\alpha,\beta>0$ be two real numbers
with the property:
\begin{equation}
\alpha D(P_1)=\beta D(P_2).
\label{pi=typeab}
\end{equation}
If, for each $k=1,\dots,4$, a mediator $M_k$ in $P_k\mathcal{A}P_k$ is given,
then the system
$\Gamma=(A,B;V,W)$, defined by
\begin{align*}
A&=\alpha P_1+\beta M_1+\alpha M_4;& B&=-\beta M_1-\alpha M_4;\\
V&=\beta P_2+\alpha M_2+\beta M_3;& W&=-\beta M_3-\alpha M_2.
\end{align*}
is a $2$-folding of $\alpha P_1$ as $\beta P_2$.
\end{lemma}
\begin{proof}
Consider the numbers $a=D(P_2)$ and $b=D(P_1)$, and denote the common value
$a /\alpha = b /\beta$ by $\lambda$.
Since, for every $j\in\{1,\dots,4\}$ we have
$$
q_{P_j\mathcal{A}P_j}(X)=\frac{q_{\mathcal{A}}(X)}{D(P_j)},
\,\,\,\forall\,X\in P_j\mathcal{A}P_j,
$$
we get the equalities
\begin{align}
q_{\mathcal{A}}(M_1^k)&=
q_{\mathcal{A}}(M_3^k)=\frac{b}{k+1},\label{qMX}\\
q_{\mathcal{A}}(M_2^k)&=
q_{\mathcal{A}}(M_4^k)=\frac{a}{k+1},\label{qNY}
\end{align}
for all integers $k\geq 1$.

Note that, since $P_1,\dots,P_4,M_1,\dots,M_4$ commute, the
elements $A,B,V,W$ also commute.

Using the obvious orthogonality relations
$(\alpha P_1+\beta M_1)\perp (\alpha M_4)$,
$(\beta M_1)\perp(\alpha M_4)$,
$(\beta P_2+\alpha M_2)\perp(\beta M_3)$, and $(\alpha M_2)\perp(\beta M_3)$, one has the
the equalities:
\begin{align}
A^k&=\alpha^k M_4^k+(\alpha P_1+\beta M_1)^k=\alpha^k\big[M_4^k+P_1\big]
+\sum_{j=1}^k\binom{k}{j}\alpha^{k-j}\beta^j
M_1^j,
\label{ak}\\
B^k&=(-1)^k\big[\beta^kM_1^k+\alpha^kM_4^k\big],
\label{bk}\\
V^k&=\beta^k M_3^k+(\beta P_2+\alpha M_2)^k=\beta^k\big[M_3^k+ P_2\big]+
\sum_{j=1}^k\binom{k}{j}\beta^{k-j}\alpha^jM_2^j,
\label{vk}\\
W^k&=(-1)^k\big[\alpha^kM_2^k+\beta^kM_3^k\big],\label{wk}
\end{align}
for every integer $k\geq 1$. It is also pretty obvious that $A,B\perp V,W$, and
we have the equalities $A+B=\alpha P_1$ and $V+W=\beta P_2$,
so in order to finish the
proof, we are left to show that $A\sim V$, and $B\sim W$. For this purpose we
use Theorem 2.1, which means that it suffices to prove the equalities
$q_{\mathcal{A}}(A^k)=q_{\mathcal{A}}(V^k)$ and
$q_{\mathcal{A}}(B^k)=q_{\mathcal{A}}(W^k)$, $\forall\,k\in\mathbb{N}$.
These equalities are proven by direct computation, as follows.

For $B$ and $W$ the equality follows from \eqref{qMX} and \eqref{qNY}, which
immediately give:
$$
q_{\mathcal{A}}(B^k)=q_{\mathcal{A}}(W^k)=
\frac{(-1)^k[\alpha^k a+\beta^k b]}{k+1}.$$
For $A$ and $V$,
again using \eqref{qMX} and \eqref{qNY},
we have
\begin{align*}
q_{\mathcal{A}}(A^k)&=\alpha^k \big[q_{\mathcal{A}}(M_4^k)+b\big]
+\sum_{j=1}^k\binom{k}{j}\alpha^{k-j}\beta^jq_{\mathcal{A}}(M_1^j)=\\
&=\frac{\alpha^k a}{k+1}+b\bigg[\alpha^k+\sum_{j=1}^k\binom{k}{j}\frac{
\alpha^{k-j}\beta^{j}}{j+1}\bigg]=
\frac{\alpha^k a}{k+1}+b\bigg[\sum_{j=0}^k\binom{k}{j}\frac{
\alpha^{k-j}\beta^j}{j+1}\bigg];\\
q_{\mathcal{A}}(V^k)&=\beta^k \big[q_{\mathcal{A}}(M_3^k)+a\big]
+\sum_{j=1}^k\binom{k}{j}\beta^{k-j}\alpha^jq_{\mathcal{A}}(M_2^j)
=\\
&=\frac{\beta^k b}{k+1}+a\bigg[\beta^k+\sum_{j=1}^k\binom{k}{j}\frac{
\beta^{k-j}\alpha^{j}}{j+1}\bigg]=
\frac{\beta^k b}{k+1}+a\bigg[\sum_{j=0}^k\binom{k}{j}\frac{
\beta^{k-j}\alpha^j}{j+1}\bigg].
\end{align*}
Replacing $a=\lambda\alpha$ and $b=\beta\lambda$, the above computations continue
as
\begin{align*}
q_{\mathcal{A}}(A^k)&=
\lambda\bigg(\frac{\alpha^{k+1}}{k+1}+\sum_{j=0}^k\binom{k}{j}\frac{
\alpha^{k-j}\beta^{j+1}}{j+1}\bigg)=\\
&=\lambda\bigg(\frac{\alpha^{k+1}}{k+1}+\int_0^\beta\big[\sum_{j=0}^k
\binom{k}{j}\alpha^{k-j}t^j\big]\,dt\bigg)=\\
&=\lambda\bigg(\frac{\alpha^{k+1}}{k+1}+\int_0^\beta(\alpha+t)^k\,dt\bigg)=
\frac{\lambda(\alpha+\beta)^{k+1}}{k+1},\\
q_{\mathcal{A}}(V^k)&=
\lambda\bigg(\frac{\beta^{k+1}}{k+1}+\sum_{j=0}^k\binom{k}{j}\frac{
\beta^{k-j}\alpha^{j+1}}{j+1}\bigg)=\\
&=\lambda\bigg(\frac{\beta^{k+1}}{k+1}+\int_0^\alpha\big[\sum_{j=0}^k
\binom{k}{j}\beta^{k-j}t^j\big]\,dt\bigg)=\\
&=\lambda\bigg(\frac{\beta^{k+1}}{k+1}+\int_0^\alpha(\beta+t)^k\,dt\bigg)=
\frac{\lambda(\alpha+\beta)^{k+1}}{k+1},
\end{align*}
so we indeed have the equality $q_{\mathcal{A}}(A^k)=q_{\mathcal{A}}(V^k)$.
\end{proof}

\begin{convention}
For the remainder of this section we are going to work under the 
following assumptions:
{\em We fix $\mathcal{A}$ to be an AW*-factor of type \twone. We fix a thick
AW*-subfactor $\mathcal{B}$ of type \twone (which exists by Example 3.2). We fix
a $\mathcal{B}$-mediator $M$ in $\mathcal{A}$.}
\end{convention}

\begin{notation}
Let $\pi=(P,Q;P',Q')\in\boldsymbol{\Pi}(\mathcal{B})$, and let
$\alpha$, $\beta$ be positive real numbers. We define the system
$\Gamma^{\alpha\beta}(\pi)=(A,B;V,W)\subset\mathcal{A}_{sa}$ by:
\begin{align*}
A&=P(\alpha I+\beta M)+\alpha Q'M;& B&=-\beta PM-\alpha Q'M;\\
V&=Q(\beta I+\alpha M)+\beta P'M;& W&=-\beta P'M-\alpha QM.
\end{align*}
\end{notation}

\begin{proposition}
If $\pi=(P,Q;P',Q')\in\boldsymbol{\Pi}(\mathcal{B})$, and if the real numbers
$\alpha,\beta>0$ satisfy the condition:
\begin{equation}
\alpha D(P)=\beta D(Q),
\label{absup}
\end{equation}
then $\Gamma^{\alpha\beta}(\pi)$ is a $2$-folding of
$\alpha P$ as $\beta Q$.
\end{proposition}
\begin{proof}
By Lemma 4.1, all we must show is the fact that
$PM$ is a mediator in $P\mathcal{A}P$,
$P'M$ is a mediator in $P'\mathcal{A}P'$,
$QM$ is a mediator in $Q\mathcal{A}Q$,
and $Q'M$ is a mediator in $Q'\mathcal{A}Q$.
But this is obvious, since $P,P',Q,Q'$ all belong to $\mathcal{B}$,
and $M$ is a $\mathcal{B}$-mediator.
\end{proof}

For the purpose of a smooth exposition, we isolate the hypothesis of
the above result as follows.
 
\begin{definition}
Given a superprojection $\pi=(P,Q;P',Q')\in\boldsymbol{\Pi}(\mathcal{B})$, and
two real numbers $\alpha,\beta>0$, we declare $\pi$ to be of
{\em of type $\alpha|\beta$} -- or say
$\pi$ is an {\em $\alpha|\beta$-superprojection} -- if $\pi$ satisfies condition
\eqref{absup}.
(The reason we use the notation $\alpha|\beta$ is the fact that the feature we
are interested in does not change if both $\alpha$ and $\beta$ are multiplied by
a factor.)
\dontwrite{
We denote by
$\boldsymbol{\Pi}^{\alpha|\beta}(\mathcal{B})$ the set of all
$\alpha|\beta$-superprojections in $\mathcal{B}$.
}
\end{definition}

\begin{theorem}[Local Folding]
Let $P,Q\in\mathbf{P}(\mathcal{B})$ be two projections with $Q\leq P$, let
$X\in P\mathcal{B}P$, be a positive element, and let $\beta>0$ be a real number
with the following properties
\begin{itemize}
\item[(i)] $X\perp Q$;
\item[(ii)] $q_{\mathcal{A}}(X)=\beta D(Q)$;
\item[(iii)] $D(P)\geq 2\big[D\big(\mathbf{s}(X)\big)+D(Q)\big]$.
\end{itemize}
Then there exists a $2$-folding $\Phi$, of $X$ as $\beta Q$, with
$\mathbf{s}(\Phi)\leq P$.
\end{theorem}
\begin{proof}
We begin the proof by fixing some notations.

Denote for simplicity $\mathbf{s}(X)$ by $S$, and $D(S)$ by $\delta$.
Use Proposition 3.6 to find a scale $\mathcal{E}=(E,[0,\delta])$ in
$\mathcal{B}$ with $E(\delta)=S$, and a non-decreasing function
$f:[0,\delta]\to \big[0,\|X\|\big]$ such that $X=\int_0^\delta f\,d\mathcal{E}$.

Since $S\perp Q$, we know that $S+Q\leq P$, and moreover
$D(P)\geq 2\big[D(S)+D(Q)\big]$. In particular, there exist two more
projections $S',Q'\in\mathbf{P}(\mathcal{B})$, with $S',Q'\leq P$, such that
$\sigma=(S,Q;S',Q')$ is a superprojection. Let then $\mathcal{E}'=(E',[0,\delta])$ be a
scale in $\mathcal{B}$, with $E'(\delta)=S'$, and define the
element $X'=\int_0^\delta f\,d\mathcal{E}'$.

Let us define, for any closed subinterval $J=[a,b]$ of $[0,\delta]$, the number
$\alpha_J=\inf_{t\in J}f(t)$ (note that $0\leq \alpha_J\leq\|X\|$), and the projections
$E_J=E(b)-E(a)$, and $E_J'=E'(b)-E'(a)$. (Of course, both
$E_J$ and $E_J'$ belong to $\mathcal{B}$, they are orthogonal -- since
$E_J\leq R$ and $E_J'\leq R'$, and they are equivalent, since
$D(E_J)=D(E_J')=b-a$.)

We also fix a sequence of partitions
$(\Delta_n)_{n=1}^\infty\subset\mathfrak{P}[0,\delta]$, such that
\begin{itemize}
\item[(a)] $\Delta_1\subset\Delta_2\subset\dots$;
\item[(b)] $\lim_{n\to\infty}L(f,\Delta_n)=\int_0^\infty f(t)\,dt$.
\end{itemize}
In fact, we can also assume that $\Delta_1=[0<\delta]$, and
\begin{itemize}
\item[(c)] {\em for every $n\geq 1$, the partition $\Delta_{n+1}$ is obtained
by subdividing exactly one interval in $\Delta_n$ into two sub-intervals.}
\end{itemize}
In other words, if $\Delta_n=[0=t_0<t_1<\dots< t_n=\delta]$, then
$\Delta_{n+1}=[0=s_0<s_1\dots<s_{n+1}=\delta]$, with
$\{t_0,t_1,\dots,t_n\}\subset\{s_0,s_1,\dots,s_{n+1}\}$. (The fact that
$\Delta_n$ consists of a partition into $n$ intervals is no coincidence.)

For every $n\geq 1$, let $\mathcal{J}_n$ be the set of intervals determined by $\Delta_n$.
(Namely, if $\Delta_n=[0=t_0<t_1<\dots <t_n=\delta]$, then
$\mathcal{J}_n=\big\{[t_{i-1},t_i]\,:\,i=1,\dots,n\big\}$.)
With this notation, $\mathcal{J}_{n+1}$ is obtained from $\mathcal{J}_n$ by
splitting (exactly) one of its intervals -- denoted $J_n$ -- into two sub-intervals,
denoted $L_n$ (the left one) and $R_n$ (the right one), so if, say $J_n=[a,b]$, then
$L_n=[a,c]$ and $R_n=[c,b]$ for some $a<c<b$. With this notation, we have:
$\mathcal{J}_{n+1}=(\mathcal{J}_n\smallsetminus\{J_n\})\cup\{L_n,R_n\}$.

Denote for simplicity $L_{\mathcal{E}}(f,\Delta_n)$ by $X_n$, and
$L_{\mathcal{E}'}(f,\Delta_n)$ by $X'_n$.
With these notations, one obviously has the equalities
\begin{equation}
X_n=\sum_{J\in\mathcal{J}_n}\alpha_J E_J\text{ and }
X'_n=\sum_{J\in\mathcal{J}_n}\alpha_J E'_J,\,\,\,\forall\,
n\in\mathbb{N},
\label{Xn=}
\end{equation}
with $0\leq X_n\leq X$ and $0\leq X'_n\leq X'$.

\begin{claimone}
Let $\mathcal{J}=\bigcup_{n=1}^\infty \mathcal{J}_n$.
There exist two maps
$\mathcal{J}\ni J\longmapsto Q_J\in\mathbf{P}(\mathcal{B})$ and
$\mathcal{J}\ni J\longmapsto Q'_J\in\mathbf{P}(\mathcal{B})$, with the following
properties:
\begin{itemize}
\item[{\sc (a)}] $Q_J\leq Q$ and $Q'_J\leq Q'$ (hence $Q_J\perp Q'_J$),
$\forall\,J\in\mathcal{J}$;
\item[{\sc (b)}] $Q_J\sim Q'_J$, and $\alpha_J D(E_J)=\beta D(Q_J)$,
$\forall\,J\in\mathcal{J}$;
\item[{\sc (c)}] if $J,K\in\mathcal{J}$ are essentially disjoint
(i.e. $J\cap K$ has at most one point), then $Q_J\perp Q_K$ and
$Q'_J\perp Q'_K$;
\item[{\sc (d)}] $Q_{J_n}\geq Q_{L_n}$ and $Q'_{J_n}\geq Q'_{L_n}$,
$\forall\,n\geq 1$;
\item[{\sc (e)}] $Q_{L_n}+Q_{R_n}\geq Q_{J_n}$ and
$Q'_{L_n}+Q'_{R_n}\geq Q'_{J_n}$,
$\forall\,n\geq 1$.
\end{itemize}
\end{claimone}
The two maps will be defined recursively. Put $\tilde{\mathcal{J}}_k=
\mathcal{J}_{k-1}\cup\mathcal{J}_k$ (with the convention $\mathcal{J}_0=\varnothing$),
so that we still have $\mathcal{J}=\bigcup_{k=1}^\infty \tilde{\mathcal{J}}_k$, and
$\tilde{\mathcal{J}}_{k+1}=\tilde{\mathcal{J}}_k\cup\{L_k,R_k\}$
(disjoint union). Of course, $\tilde{\mathcal{J}}_1=\mathcal{J}_1=\{J_1\}$, where
$J_1=[0,\delta]$. Start off by defining $Q_{J_1}$ and $Q'_{J_1}$ to be arbitrary projections
with $Q_{J_1}\leq Q$ and $Q'_{J_1}\leq Q'$, such that 
$$D(Q_{J_1})=D(Q'_{J_1})=\frac{\alpha_{J_1}D(E_{J_1})}{\beta}.$$
This is possible, since
$$\frac{\alpha_{J_1}D(E_{J_1})}{\beta}=\frac{q_{\mathcal{A}}(X_1)}{\beta}
\leq \frac{q_{\mathcal{A}}(X)}{\beta}=D(Q)=D(Q').$$
Assume now the projections $Q_J$ and $Q'_J$ are defined for all
$J\in \tilde{\mathcal{J}}_k$, and they satisfy conditions
{\sc (a)(b)(c)} (with $\tilde{\mathcal{J}}_k$ in place of $\mathcal{J}$),
and conditions {\sc (d)} and {\sc (e)} for all
$n<k$, and let us indicate how to construct the ``new'' projections
$Q_{L_k}$, $Q'_{L_k}$, $Q_{R_k}$, and $Q'_{R_k}$.

Define the elements $F_k=\sum_{J\in\mathcal{J}_k}Q_J$ and
$F'_k=\sum_{J\in\mathcal{J}_k}Q'_J$ in $\mathcal{B}$.
First of all, using {\sc (b) (c)}, it follows that $F_k$ and $F'_k$
 are projections in $\mathcal{B}$, with
$F_k\leq Q$, $F'_k\leq Q'$, and
\begin{equation}
\beta D(F_k)=\beta D(F'_k)=\sum_{J\in\mathcal{J}_k}\alpha_JD(E_J)=q_{\mathcal{A}}(X_k).
\label{DF=Xk}
\end{equation}
Secondly,
since $f$ is non-decreasing, we have:
\begin{equation}
\alpha_{L_k}=\alpha_{J_k}\text{ and }
\alpha_{R_k}\geq \alpha_{J_k}.
\label{alphaJ-f}
\end{equation}
Using {\sc (b)} for $J=J_k$,
it follows that
$$\frac{\alpha_{L_k}D(E_{L_k})}{\beta}=
\frac{\alpha_{J_k}D(E_{L_k})}{\beta}\leq
\frac{\alpha_{J_k}D(E_{J_k})}{\beta}=D(Q_{J_k})=D(Q'_{J_k}),$$
so we can choose two projections $Q_{L_k},Q'_{L_k}
\in\mathbf{P}(\mathcal{B})$, with
$Q_{L_k}\leq Q_{J_k}$, $Q'_{L_k}\leq Q'_{J_k}$, and
\begin{equation}
D(Q_{L_k})=D(Q'_{L_k})=\frac{\alpha_{L_k}D(E_{L_k})}{\beta}.
\label{DQr}\end{equation}
Let us observe that, using \eqref{DF=Xk} and the obvious equality
$$X_{k+1}-X_k=[\alpha_{R_k}-\alpha_{J_k}]\cdot E_{R_k},$$
one has (all terms commute):
\begin{gather*}
[\alpha_{R_k}-\alpha_{J_k}]\cdot D(E_{R_k})=
q_{\mathcal{A}}(X_{k+1}-X_k)=
q_{\mathcal{A}}(X_{k+1})-
q_{\mathcal{A}}(X_k)\leq \\
\leq q_{\mathcal{A}}(X)-
q_{\mathcal{A}}(X_k)=\beta D(Q)-\beta D(F_k)=\beta D(Q-F_k).
\end{gather*}
In particular, one has the inequalities
$$D(Q-F_k)=D(Q'-F'_k)\geq\frac{[\alpha_{R_k}-\alpha_{J_k}]
\cdot D(E_{R_k})}{\beta},$$
so there exist projections $G,G'\in\mathbf{P}(\mathcal{B})$, with
$G\leq Q-F_k$, $G'\leq Q'-F'_k$, such that
$$D(G)=D(G')=\frac{[\alpha_{R_k}-\alpha_{J_k}]\cdot D(E_{R_k})}{\beta}.$$
By construction, we have
\begin{equation}
G\perp F_k\text{ and }G'\perp F'_k.
\label{GperpF}
\end{equation}
Since $0\leq Q_{J_k}-Q_{L_k}\leq
Q_{J_k}\leq F_k$ and
$0\leq Q'_{J_k}-Q'_{L_k}\leq
Q'_{J_k}\leq F'_k$, using \eqref{GperpF}, we have $G\perp (Q_{J_k}-Q_{L_k})$
and
$G'\perp (Q'_{J_k}-Q'_{L_k})$. We can then define then the
projections $Q_{R_k}=G+(Q_{J_k}-Q_{L_k})$ and
$Q'_{R_k}=G'+(Q'_{J_k}-Q'_{L_k})$.

We now check conditions {\sc (a)(b)(c)} with $\tilde{\mathcal{J}}_{k+1}$ in
place of $\mathcal{J}$, and conditions {\sc (d)(e)} with $n=k$.

To check condition {\sc (a)}, we only need to consider the ``new'' intervals,
namely the cases $J=L_k,R_k$, which are obvious by construction. 

To check condition {\sc (b)} we need to prove the equalities
\begin{align}
\alpha_{L_k} D(E_{L_k})&=\beta D(Q_{L_k})=
\beta D(Q'_{L_k}),\label{b-for-l}\\
\alpha_{R_k} D(E_{R_k})&=\beta D(Q_{R_k})=
\beta D(Q'_{R_k}).\label{b-for-r}
\end{align}
The equalities \eqref{b-for-l} are trivial using \eqref{DQr}.
To prove the equalities \eqref{b-for-r}, we first
notice that by construction
we have 
\begin{align*}
D(Q_{R_k})&=D(G)+D(Q_{J_k}-Q_{L_k})=D(G)+D(Q_{J_k})-D(Q_{L_k})=\\
&=
\frac{[\alpha_{R_k}-\alpha_{J_k}]\cdot D(E_{R_k})+\alpha_{J_k}D(E_{J_k})-
\alpha_{L_k}D(E_{L_k})}{\beta}.
\end{align*}
Using the equalities $\alpha_{L_k}=\alpha_{J_k}$, as well as
$D(E_{J_k})=D(E_{L_k})+D(E_{R_k})$, the above computation
continues as
\begin{align*}
D(Q_{R_k})&=
\frac{[\alpha_{R_k}-\alpha_{J_k}]\cdot D(E_{R_k})+\alpha_{J_k}D(E_{J_k})-
\alpha_{J_k}\cdot\big[D(E_{J_k})-D(E_{L_k})\big]}{\beta}=\\
&=\frac{\alpha_{R_k} D(E_{R_k})}{\beta}=\dots
\text{ (same computation for $Q'_{R_k}$) }\hdots=
D(Q'_{R_k}),
\end{align*}
thus proving \eqref{b-for-r}.

To prove condition {\sc (c)} we only need to examine the ``new'' cases, which are
\begin{align}
Q_{L_k}\perp Q_{R_k}&\text{, and }Q'_{L_k}\perp Q'_{R_k};
\label{rperpl}\\
Q_{L_k}\perp Q_J&\text{, and }
Q'_{L_k}\perp Q'_J,\,\,\,\forall\,J\in\mathcal{J}_k\smallsetminus\{J_k\}.
\label{lperpJ};\\
Q_{R_k}\perp Q_J&\text{, and }
Q'_{R_k}\perp Q'_J,\,\,\,\forall\,J\in\mathcal{J}_k\smallsetminus\{J_k\}.
\label{rperpJ}.
\end{align}
The orthogonality relations \eqref{rperpl} follow from \eqref{GperpF}, which together
with the obvious inequalities $Q_{L_k}\leq Q_{J_k}\leq F_k$ and
$Q'_{L_k}\leq Q'_{J_k}\leq F'_k$ force $G\perp Q_{L_k}$ and
$G'\perp Q'_{L_k}$.
To prove \eqref{lperpJ} and \eqref{rperpJ} we simply observe that
$$Q_{L_k}+Q_{R_k}=Q_{J_k}+G\text{ and }
Q'_{L_k}+Q'_{R_k}=Q_{J_k}+G,$$
and by \eqref{GperpF} we also have
$$G\perp Q_J\text{ and }G'\perp Q'_J,\,\,\,\forall\,J\in\mathcal{J}_k.$$

Conditions {\sc (d)(e)} -- for $k=n$ -- are automatically satisfied, by construction.

Having proven Claim 1, we continue the proof of the Theorem,
by fixing the maps
$(Q_J)_{J\in\mathcal{J}}$ and $(Q'_J)_{J\in\mathcal{J}}$ as above.
It is obvious that, for every $J\in\mathcal{J}$, the system
$\pi_J=(E_J,Q_J;E'_J,Q'_J)$ is an $\alpha_J|\beta$-superprojection.
(By construction $\pi_J\leq\sigma$.)
Use Proposition 4.1 to define, for each $J\in\mathcal{J}$, a $2$-folding
$\Gamma_J=(A^J,B^J;V^J,W^J)$,
by
\begin{align*}
A^J&=E_J(\alpha_J I+\beta M)+\alpha_J Q'_JM;&
B^J&=-\beta E_JM-\alpha_J Q'_JM;\\
V^J&=Q_J(\beta I+\alpha_J M)+\beta E'_JM;&
W^J&=-\beta E'_JM-\alpha_J Q_JM.
\end{align*}
\begin{claimtwo}
The $2$-foldings $(\Gamma_J)_{J\in\mathcal{J}}$ have the following properties:
\begin{itemize}
\item $\|\Gamma_J\|\leq 2\big[\beta+\|X\|\big]$, $\forall\,J\in\mathcal{J}$.
\item if $J,K\in\mathcal{J}$ are essentially disjoint, then
$\mathbf{s}(\Gamma_J)\perp\mathbf{s}(\Gamma_K)$;
\end{itemize}
\end{claimtwo}
The first assertion is trivial, since $0\leq\alpha_J\leq\|X\|$, $\forall\,J\in\mathcal{J}$.

To prove the second property, we start with two intervals $J,K\in\mathcal{J}$ that
are essentially disjoint, and we notice that the collection
$\{E_J,E'_J,E_K,E'_K,Q_J,Q'_J,Q_K,Q'_K\}$ is orthogonal. (This can be done
by ``groupping,'' observing that, since we have
$E_J,E_K\leq S$, $E_J,E_K\leq S'$,
$Q_J,Q_K\leq Q$, and $Q'_J,Q'_K\leq Q$, with
$\{S,S',Q,Q'\}$ orthogonal, all we must show are the orthogonality relations:
$E_J\perp E_K$, $E'_J\perp E_K$, $Q_J\perp Q_K$, and $Q'_J\perp Q'_K$, which
are obvious.) The first stament is then clear, since we have the
following (obvious) inequalities:
\begin{align*}
\mathbf{s}(\Gamma^J)&\leq E_J+E'_J+Q_J+Q'_J,\\
\mathbf{s}(\Gamma^K)&\leq E_K+E'_K+Q_K+Q'_K,
\end{align*}
with $(E_J+E'_J+Q_J+Q'_J)\perp (E_K+E'_K+Q_K+Q'_K)$.

Having proven Claim 2, let us define now, for every integer $n\geq 1$, the
system
$\Phi_n=(A_n,B_n;V_n,W_n)$, where
$A_n=\sum_{J\in\mathcal{J}_n}A^J$,
$B_n=\sum_{J\in\mathcal{J}_n}B^J$,
$V_n=\sum_{J\in\mathcal{J}_n}V^J$, and
$W_n=\sum_{J\in\mathcal{J}_n}W^J$.

\begin{claimthree}
For every $n\in\mathbb{N}$, the system $\Phi_n$ is a
$2$-folding of $X_n$ as $\beta F_n$, where
$F_n=\sum_{J\in \mathcal{J}_n}Q_J$. Moreover, one has
$\|\Phi_n\|\leq 2\big[\beta+\|X\|\big]$, $\forall\,n\in\mathbb{N}$.
\end{claimthree}
This follows from Remark 4.1, combined with Claim 2, and
the fact that all the intervals in $\mathcal{J}_n$ are
essentially disjoint. By construction, we have
\begin{align*}
A_n+B_n&=\sum_{J\in\mathcal{J}_n}\alpha_JE_J=X_n,\\
V_n+W_n&=\beta\sum_{J\in\mathcal{J}_n}Q_J=\beta F_n,
\end{align*}
so $\Phi_n$ is indeed a $2$-folding of $X_n$ as $\beta F_n$.

\begin{claimfour}
The sequences $(A_n)_{n\in\mathbb{N}}$,
$(B_n)_{n\in\mathbb{N}}$, $(V_n)_{n\in\mathbb{N}}$, and
$(W_n)_{n\in\mathbb{N}}$ have the following properties:
\begin{itemize}
\item they are all bounded;
\item they are jointly abelian;
\item they all all lie in $(P\mathcal{A}P)_{sa}$; more precisely, for every
$n\in\mathbb{N}$, one has the
inequalities:
\begin{align}
\mathbf{s}(A_n)\leq S+Q',&\quad\mathbf{s}(B_n)\leq S+Q',\label{sAB}\\
\mathbf{s}(V_n)\leq S'+Q,&\quad\mathbf{s}(W_n)\leq S'+Q.
\label{sVW}
\end{align} 
\end{itemize}
Moreover:
\begin{itemize}
\item the sequences $(A_n)_{n\in\mathbb{N}}$ and $(V_n)_{n\in\mathbb{N}}$
are non-decreasing;
\item the sequences $(B_n)_{n\in\mathbb{N}}$ and $(W_n)_{n\in\mathbb{N}}$
are non-increasing.
\end{itemize}
\end{claimfour}
The first assertion is quite clear.
The inequalities \eqref{sAB} \eqref{sVW} are also clear, and they imply the
third assertion.
The fact that the four sequences
are jointly abelian follows from the fact that the collection
\begin{equation}
\mathcal{C}=\mathcal{E}\cup\mathcal{E}'\cup\{M\}\cup\{Q_J\,:\,J\in\mathcal{J}\}\cup
\{Q'_J\,:\,J\in\mathcal{J}\}\subset (P\mathcal{A}P)_{sa}\label{colC}
\end{equation}
is abelian, and obviously all sequences lie in the abelian von Neumann
subalgebra $\mathcal{N}=\mathcal{C}''$.

We proceed now with the proof of the monotonicity features. Fix some
$k\in\mathbb{N}$, and let us compare $A_{k+1}$ with $A_k$,
$B_{k+1}$ with $B_n$, $V_{k+1}$ with $V_k$, and $W_{k+1}$ with $W_k$.
Since
$\mathcal{J}_k=(\mathcal{J}_n\smallsetminus\{J_k\})\cup\{L_k,R_k\}$, we have
\begin{align*}
&A_{k+1}-A_k=A^{L_k}+A^{R_k}-A^{J_k}=
\big[E_{L_k}(\alpha_{L_k}I+\beta M)+
\alpha_{L_k} Q'_{L_k}M\big]+\\
&\qquad
+
\big[E_{R_k}(\alpha_{R_k} I+\beta M)+
\alpha_{R_k} Q'_{R_k}M\big]-\big[E_{J_k} (\alpha_{J_k}I+\beta M)+
\alpha_{J_k} Q'_{J_n}M\big];\\
&B_k-B_{k+1}=B^{J_k}-B^{L_k}-B^{R_k}=
\big[\beta E'_{L_k}M+\alpha_{L_k}Q'_{L_k}M\big]+\\
&\qquad
+\big[\beta E'_{R_k}M+\alpha_{R_k}Q'_{R_k}M\big]
-\big[\beta E'_{J_k}M+\alpha_{J_k}Q'_{J_k}M\big];\\
&V_{k+1}-V_k=V^{L_k}+V^{R_k}-V^{J_k}=
\big[Q_{L_k}(\beta I+\alpha_{L_k}M)+\beta
E'_{L_k}M\big]+\\
&\qquad
+\big[Q_{R_k}(\beta I+\alpha_{R_k}M)+\beta
E'_{R_k}M\big]-
\big[Q_{J_k}(\beta I+\alpha_{J_k}M)+\beta
E'_{J_k}M\big];\\
&W_k-W_{k+1}=W^{J_k}-W^{L_k}-W^{R_k}=
\big[\beta E'_{L_k}M+\alpha_{L_k}Q_{L_k}M\big]
+\\
&\qquad
+\big[\beta E'_{R_k}M+\alpha_{R_k}Q_{R_k}M\big]-
\big[\beta E'_{J_k}M+\alpha_{J_k}Q_{J_k}M\big].
\end{align*}
Using \eqref{alphaJ-f} and property {\sc (e)} from Claim 1,
combined with the equalities
$E_{J_k}=E_{J_k^\ell}+E_{J_k^r}$ and $E'_{J_k}=E'_{J_k^\ell}+E'_{J_k^r}$,
we can continue with:
\begin{align*}
&A_{k+1}-A_k=[\alpha_{R_k}-\alpha_{J_k}]E_{R_k}+
[\alpha_{J_k}Q_{L_k}+\alpha_{R_k}Q_{R_k}-\alpha_{J_k}Q_{J_k}]M\geq\\
&\qquad\geq \alpha_{J_k}[Q_{L_k}+Q_{R_k}-Q_{J_k}]M\geq 0;\\
&B_k-B_{k+1}=[\alpha_{J_k}Q'_{L_k}+\alpha_{R_k}Q'_{R_k}-\alpha_{J_k}Q'_{J_k}]M
\geq\alpha_{J_k}[Q'_{L_k}+Q'_{R_k}-Q'_{J_k}]M\geq 0;\\
&V_{k+1}-V_k= \beta[Q_{L_k}+Q_{R_k}-Q_{J_k}]+[\alpha_{J_k}Q_{L_k}+
\alpha_{R_k}Q_{R_k}-\alpha_{J_k}Q_{J_k}]M\geq \\
&\qquad\geq \alpha_{J_k}[Q_{L_k}+Q_{R_k}-Q_{J_k}]M\geq 0;\\
&W_k-W_{k+1}=[\alpha_{J_k}Q_{L_k}+\alpha_{R_k}Q_{R_k}-\alpha_{J_k}Q_{J_k}]M\geq
\alpha_{J_k}[Q_{L_k}+Q_{R_k}-Q_{J_k}]M\geq 0.
\end{align*}

Having proven Claim 4,  we now use Lemma 1.1, which gives the existence of the
weak limits $A=\wlim_{n\to\infty}A_n$, $B=\wlim_{n\to\infty}B_n$,
$V=\wlim_{n\to\infty}V_n$, and $W=\wlim_{n\to\infty}W_n$.
The proof of the Theorem will be finished, once we prove the following.
\begin{claimfive}
The system $\Phi=(A,B;V,W)$ is a $2$-folding of $X$ as $\beta Q$, with
$\mathbf{s}(\Phi)\leq P$.
\end{claimfive}
First of all, since all four sequences lie in the abelian von Neumann
algebra $\mathcal{N}=\mathcal{C}''$, with $\mathcal{C}$ defined by
\eqref{colC}, it follows that $A$, $B$, $V$, and $W$ all belong to
$\mathcal{N}$. In particular, these four self-adjoint elements commute.

Note also that, since all sequences lie in $P\mathcal{A}P$, it follows that
$A,B,V,W\in P\mathcal{A}P$.

Using Proposition 2.3, in conjunction with Claim 3, it is obvious that
$A\sim V$ and $B\sim W$.

Moreover, by Remark 1.4 and by Claim 3, we also have the equalities
\begin{align}
A+B&=\wlim_{n\to\infty}(A_n+B_n)=\wlim_{n\to\infty}X_n=X;
\label{A+B=X}\\
V+W&=\wlim_{n\to\infty}(V_n+W_n)=\beta\wlim_{n\to\infty}F_n.
\label{V+W=F}
\end{align}
Since
$$F_{k+1}-F_k=Q^{L_k}+Q^{R_k}-Q_{J_k}\geq 0,\,\,\,\forall\,k\in\mathbb{N},$$
we know that $\wlim_{n\to\infty}F_n$ is in fact a projection. Moreover, since
we also know that $F_n\leq Q$, i.e. $F_n\in Q\mathcal{A}Q$,
$\forall\,n\in\mathbb{N}$, it follows that
$F\in Q\mathcal{A}Q$, so $F\leq Q$. By \eqref{A+B=X} and \eqref{V+W=F},
 combined with Lemma 1.1, we know that
\begin{align*}
\beta D(F)&=\lim_{n\to\infty}q_{\mathcal{A}}(\beta F_n)=
\lim_{n\to\infty}\big[q_{\mathcal{A}}(V_n)+q_{\mathcal{A}}(W_n)\big]=\\
&=\lim_{n\to\infty}\big[q_{\mathcal{A}}(A_n)+q_{\mathcal{A}}(B_n)\big]=
\lim_{n\to\infty}q_{\mathcal{A}}(X_n)=
q_{\mathcal{A}}(X).\end{align*}
This forces, of course $\beta D(F)=\beta D(Q)$, and then the condition $F\leq Q$
(combined with $\beta >0$) will force $F=Q$.

At this point the only properties left to be proven are the orthogonality relations
$A\perp V$, $A\perp W$, $B\perp V$, and $B\perp W$.
For this purpose we use \eqref{sAB} and \eqref{sVW}, to conclude that
$A,B\in (S+Q')\mathcal{A}(S+Q')$ and
$V,W\in (S'+Q)\mathcal{A}(S+Q')$, and
then everything follows from $(S+Q')\perp (S'+Q)$.
\end{proof}

\


\section{Self-adjoint elements with zero quasitrace}

In this section we prove the main results of this paper.

\dontwrite{
For technical purposes,
we introduce the following terminology.
\begin{definition}
Given a C*-algebra $\mathcal{A}$, we call an element $S\in\mathcal{A}_{sa}$ a
{\em quasisymmetry}, if its positive and negative parts $S^+$ and $S^-$ are
equivalent projections.
\end{definition}

With this terminology, we state an important application of the results from
Section 4.
}
\begin{theorem}
Let $\mathcal{A}$ be an AW*\-factor of type \twone, and let
$X\in\mathcal{A}_{sa}$ be an element with $q_{\mathcal{A}}(X)=0$.
Assume:
\begin{itemize}
\item[(i)] $D\big(\mathbf{s}(X)\big)<\frac 12$;
\item[(ii)] there is a thick AW*-subfactor $\mathcal{B}$, of type \twone,
which contains $X$.
\end{itemize}
Then there exist a spectrally symmetric element $S\in\mathcal{B}$ and
a $2$-folding of $X$ as $S$.
\end{theorem}
\begin{proof}
We assume of course that $X\neq 0$ (otherwise we can take $X_1=X_2=S=0$).

Let $E_1=\mathbf{s}(X^+)$ and $E_2=\mathbf{s}(X^-)$, so that
$E_1,E_2\in\mathbf{P}(\mathcal{B})$,
are orthogonal, and $\mathbf{s}(X)=E_1+E_2$.
By condition (ii) we know that $D(E_1)+D(E_2)<\frac 12$.
Chose then six more projections $E_3,E_4,Q_1,Q_2,Q_3,Q_4\in\mathbf{P}(\mathcal{B})$,
such that
\begin{itemize}
\item[{\sc (a)}] all eight projections $E_1,\dots,E_4,Q_1,\dots,Q_4$ are orthogonal;
\item[{\sc (b)}] $E_1+\dots +E_4+Q_1+\dots +Q_4=I$;
\item[{\sc (c)}] $E_1\sim E_3$ and $E_2\sim E_4$;
\item[{\sc (d)}] $Q_1\sim Q_2\sim Q_3\sim Q_4$.
\end{itemize}
Denote $D(Q_1)$ simply by $\delta$, so that $D(Q_k)=\delta$, $k=1,\dots,4$.

Since we are assuming $X\neq 0$, and $q_{\mathcal{A}}(X)=0$,
it follows that
$$q_{\mathcal{A}}(X^+)=q_{\mathcal{A}}(X^-)>0.$$
Denote this common value by $\alpha$, and let $\beta=\alpha/\delta$.

Consider now the projections $P_1=E_1+E_3+Q_1+Q_2$ and $P_2=I-P_1=
E_2+E_4+Q_3+Q_4$, which satisfy:
\begin{align}
D(P_1)&=2\big[D(E_1)+D(Q_1)\big];\label{DP1}\\
D(P_2)&=2\big[D(E_2)+D(Q_3)\big];\label{DP2}
\end{align}
Remark that using the above two equalities we have the following situations.
\begin{itemize}
\item[{\sc (i)}] Since $\mathbf{s}(X^+)=E_1$, it follows that:
\begin{itemize}
\item[$\bullet$] $X^+,Q_1\in P_1\mathcal{B}P_1$;
\item[$\bullet$] $\mathbf{s}(X^+)\perp Q_1$
\item[$\bullet$] $D(P_1)=2\big[D(\mathbf{s}(X^+))+D(Q_1)\big]$;
\item[$\bullet$] $q_{\mathcal{A}}(X^+)=\beta D(Q_1)$.
\end{itemize}
\item[{\sc (ii)}] Since $\mathbf{s}(X^-)=E_2$, it follows that:
\begin{itemize}
\item[$\bullet$] $X^-,Q_3\in P_2\mathcal{B}P_2$;
\item[$\bullet$] $\mathbf{s}(X^-)\perp Q_3$
\item[$\bullet$] $D(P_2)=2\big[D(\mathbf{s}(X^-))+D(Q_3)\big]$;
\item[$\bullet$] $q_{\mathcal{A}}(X^-)=\beta D(Q_3)$.
\end{itemize}
\end{itemize}
We now use Theorem 4.1 to find
\begin{itemize}
\item[{\sc (i)}] a $2$-folding $\Phi_1=(A_1,B_1;V_1,W_1)$ of $X^+$ as $\beta
Q_1$, with $\mathbf{s}(\Phi_1)\leq P_1$, and
\item[{\sc (ii)}] a $2$-folding $\Gamma=(A_2,B_2;V_2,W_2)$ of $X^-$ as $\beta
Q_2$, with $\mathbf{s}(\Gamma)\leq P_2$.
\end{itemize}
It is trivial to see that the system
$\Phi_2=(-A_2,-B_;-V_2,-W_2)$ is a $2$-folding of $-X^-$ as $-\beta Q_3$, again with
$\mathbf{s}(\Phi_2)\leq P_2$. By Remark 4.1, the system
$$\Phi_1+\Phi_2=(A_1-A_2,B_1-B_2;V_1-V_2,W_1-W_2)$$
is a $2$-folding of $X^+-X^-=X$ as $\beta Q_1-\beta Q_2$.
Obviously the element
$S=\beta Q_1-\beta Q_2$ is spectrally symmetric.
\dontwrite{
In particular (see the definition at the beginning of Section 4),
the commuting elements $X_1=(A_1-A_2)-(V_1-V_2)$ and
$X_2=(B_1-B_2)-(W_1-W_2)$ are spectrally symmetric, and they satisfy the equality
$$X_1+X_2=X-\beta Q_1+\beta Q_3,$$
so if we put $S=Q_1-Q_3$, the desired conclusion follows.
}
\end{proof}

\begin{corollary}
If $\mathcal{A}$ is an AW*-factor, and if
$X\in\mathcal{A}_{sa}$ is an element with $q_{\mathcal{A}}(X)$, satisfying the
additional conditions {\rm (i)} and {\rm (ii)} from Theorem 5.1, then
$X$ can be written as a sum of three commuting spectrally symmetric
elements in $\mathcal{A}$.
\end{corollary}
\begin{proof}
Let $S$ be a spectrally symmetric element, such that there exists
a $2$-folding $\Phi=(A_1,A_2;B_1,B_2)$ of $X$ as $S$. Then
$$X=(A_1-B_1)+(A_2-B_2)+S$$
is a sum of the desired form.
\end{proof}

The discussion from this point on is aimed at removing
condition (ii) from the hypothesis of Theorem 5.1, and relaxing condition (i)
as much as possible.

\begin{lemma}[Small Packing]
Let $\mathcal{A}$ be an AW*-factor of type \twone, let
$P\in\mathbf{P}(\mathcal{A})$ be a non-zero projection, and let $\mathcal{B}$ be an
AW*-subfactor of $P\mathcal{A}P$.
For any element $X\in\mathcal{A}_{sa}$ with $X\perp P$, and any
non-zero projection $Q\in\mathbf{P}(\mathcal{B})$, there exist five elements
$A_1,A_2,B_1,B_2,Y\in\mathcal{A}_{sa}$, with the following properties
\begin{itemize}
\item[(i)] $A_1,A_2,B_1,B_2,Y$ all commute;
\item[(ii)] $A_1\perp A_2$, $B_1\perp B_2$, $A_1\perp B_1$ and $A_2\perp B_2$;
\item[(iii)] $A_1\sim B_1$ and $A_2\sim B_2$;
\item[(iv)] $A_1,A_2,B_1\perp P$;
\item[(v)] $Y\in Q\mathcal{B}Q$;
\item[(vi)] $B_2P=PB_2=Y$;
\item[(vii)] $X=A_1+A_2-B_1-B_2+Y$.
\end{itemize}
\end{lemma}
\begin{proof}
We will assume $D_{\mathcal{A}}\big(\mathbf{s}(X)\big)>0$ (the case $X=0$ is trivial).
Let $\lambda=D_{\mathcal{A}}(P)$, so that
\begin{equation}
q_{\mathcal{A}}(B)=\lambda q_{\mathcal{B}}(B),\,\,\,\forall\,B\in\mathcal{B}.
\label{qBlambda}
\end{equation}
Let $\beta=D_{\mathcal{A}}(Q)$, so that $D_{\mathcal{B}}(Q)=\beta/\lambda$.

Fix some integer $n\geq 1$, such that
$$2n\geq\frac{D_{\mathcal{A}}\big(\mathbf{s}(X)\big)}{\beta},$$
and define the number
$$\alpha=\frac{D_{\mathcal{A}}\big(\mathbf{s}(X)\big)}{2n},$$
so that we have the equality
$D_{\mathcal{A}}\big(\mathbf{s}(X)\big)=2n\alpha$, and
$D_{\mathcal{A}}(Q)\geq \alpha$.

Using Proposition 3.5 there is a scale
$\mathcal{F}=\big(F,[0,2n\alpha]\big)$ in $\mathcal{A}$ with
$F(2n\alpha)=\mathbf{s}(X)$, and a
function $f\in\mathfrak{R}[0,2n\alpha]$, such that
$X=\int_0^{2n\alpha} f\,d\mathcal{F}$. Fix also a full
scale $\mathcal{G}=\big(G,[0,1]\big)$ in $\mathcal{B}$, that contains $Q$, so that
$D_{\mathcal{B}}\big(G(\beta/\lambda)\big)=Q$.
By construction, one has
\begin{equation}
D_{\mathcal{A}}\big(G(t/\lambda)\big)=\lambda
D_{\mathcal{B}}\big(G(t/\lambda)\big)=t,\,\,\,\forall\,t\in [0,\lambda].
\label{DAG}
\end{equation}
Since $\alpha\leq\beta$, we have $G(\alpha/\lambda)\leq Q$.

Using \eqref{DAG} it follows that the system
$\mathcal{E}=\big(E,[0,(2n+1)\alpha])$ defined by:
$$
E(t)=\left\{\begin{array}{cl}F(t)&\text{ if }0\leq t\leq 2n\alpha\\
G\big((t-2n\alpha)/\lambda\big)+F(2n\alpha)&\text{ if }
2n\alpha<t<(2n+1)\alpha\end{array}\right.
$$
is a scale in $\mathcal{A}$.
Its key features are as follows:
\begin{itemize}
\item[{\sc (I)}] $X=\int_0^{2n\alpha}f\,d\mathcal{E}$;
\item[{\sc (II)}] $E(t)\perp P$, $\forall\,t\in[0,2n\alpha]$;
\item[{\sc (III)}] $E(t)-E(2n\alpha)\in Q\mathcal{B}Q$, $\forall\,t\in[2n\alpha,(2n+1)\alpha]$.
\end{itemize}
Define the functions
$$f_k=f\big|_{[(k-1)\alpha,k\alpha]}\in\mathfrak{R}\big[(k-1)\alpha,k\alpha\big],
\,\,\,k=1,2,\dots,2n,$$
so that one has the equality
\begin{equation}
X=\sum_{k=1}^{2n}\int_{(k-1)\alpha}^{k\alpha}f_k\,d\mathcal{E}.
\label{X=sum-int}
\end{equation}
Define the sequence of functions
$g_k\in\mathfrak{R}\big[(k-1)\alpha,k\alpha\big]$, $k=1,2,\dots,2n$
starting with $g_1=f_1$, and using the recursive formula
$$g_k=f_k+\Lambda_\alpha g_{k-1},\,\,\,k=2,3,\dots,2n.$$
Here $\Lambda_\alpha:\mathfrak{R}[a,b]\to\mathfrak{R}[a+\alpha,b+\alpha]$ denotes
the translation map (see Remark 3.5).

Define now the sequences $(V_k)_{k=1}^{2n}$ and $(W_k)_{k=1}^{2n}$ by
$$V_k=\int_{(k-1)\alpha}^{k\alpha}g_k\,d\mathcal{E}
\text{ and }W_k=\int_{k\alpha}^{(k+1)\alpha}\Lambda_\alpha g_k\,d\mathcal{E},
\,\,\,k=1,2,\dots,2n.$$
The key features of these two sequences are described below.
\begin{claim}
The sequences $(V_k)_{k=1}^{2n}$ and $(W_k)_{k=1}^{2n}$
have the following properties:
\begin{itemize}
\item[{\sc (a)}] $V_1,\dots,V_{2n},W_1,\dots,W_{2n}$ all commute;
\item[{\sc (b)}] $V_i\perp V_j$ and $W_i\perp W_j$, $\forall\,i,j\in\{1,\dots,2n\}$ with
$i\neq j$;
\item[{\sc (c)}] $V_i\perp P$, $\forall\,i\in\{1,\dots,2n\}$;
\item[{\sc (d)}] $W_i\perp P$, $\forall\,i\in\{1,\dots,2n-1\}$;
\item[{\sc (e)}] $V_i\perp W_j$, $\forall\,i,j\in\{1,\dots,2n\}$, with $j\neq i-1$.
\item[{\sc (f)}] $W_{2n}\in Q\mathcal{B}Q$;
\item[{\sc (g)}] $X=\sum_{i=1}^{2n}V_i-\sum_{j=1}^{2n-1}W_j$;
\item[{\sc (a)}] $V_i\sim W_i$, $\forall\,i\in\{1,\dots,2n\}$.
\end{itemize}
\end{claim}
The first assertion is trivial.
To prove properties {\sc (a)-(e)} we define the projections
$R_k=E(k\alpha)-E((k-1)\alpha)$, $k=1,\dots, 2n+1$, and we observe that
\begin{itemize}
\item $\mathbf{s}(V_i)\leq R_i$ and $\mathbf{s}(W_i)\leq R_{i+1}$,
$\forall\,i\in\{1,\dots,2n\}$;
\item $R_i\perp R_j$, $\forall\,i,j\in\{1,\dots,2n+1\}$ with $i\neq j$;
\item $R_i\perp P$, $\forall\,i\in\{1,\dots,2n\}$;
\item $R_{2n+1}\in Q\mathcal{B}Q$
\end{itemize}

The fact that $W_{2n}=\int_{2n\alpha}^{(2n+1)\alpha}\Lambda_\alpha g_{2n}\,
d\mathcal{E}$ belongs to $Q\mathcal{B}$ follows from condition {\sc (III)} above.

Property {\sc (g)} is quite obvious from \eqref{X=sum-int}, since
$V_1=\int_0^\alpha f_1\,d\mathcal{E}$, and
$$V_k-W_{k-1}=\int_{(k-1)\alpha}^{k\alpha}
(g_k-\Lambda_\alpha g_{k-1})\,d\mathcal{E}=
\int_{(k-1)\alpha}^{k\alpha}
f_k\,d\mathcal{E},\,\,\,\forall\,k\in\{2,3,\dots,2n\}.$$

Finally, property {\sc (h)} is immediate from Corollary 3.3.

Having proven the above Claim, we now define the elements
$A_1=\sum_{k=1}^n V_{2k-1}$,
$A_2=\sum_{k=1}^n V_{2k}$,
$B_1=\sum_{j=1}^n W_{2k-1}$,
$B_2=\sum_{j=k}^nW_{2k}$, and $Y=W_{2n}$.
Using the Claim and Corollary 2.1, it is pretty obvious that
$A_1\sim B_1$ and $A_2\sim B_2$. The fact that the elements
$A_1,A_2,B_1,B_2,Y$ satisfy all the other desired conditions follows
from the Claim.
\end{proof}

\begin{theorem}
Let $\mathcal{A}$ be an AW*\-factor of type \twone, and let
$X\in\mathcal{A}_{sa}$ be an element with $q_{\mathcal{A}}(X)=0$.
If $D\big(\mathbf{s}(X)\big)< 1$,
then $X$ can be written as a sum
$X=X_1+X_2+X_3$ of three commuting
spectrally symmetric elements $X_1,X_2,X_3\in\mathcal{A}_{sa}$.
\end{theorem}
\begin{proof}
Put $P=I-\mathbf{s}(X)$, and denote the AW*-subalgebra
$P\mathcal{A}P$ by $\mathcal{A}_0$. Of course $\mathcal{A}_0$ is an AW*-factor
of type \twone. Fix some thick AW*-subfactor $\mathcal{B}$ of $\mathcal{A}_0$,
of type \twone,  as well as some projection $Q\in\mathbf{P}(\mathcal{B})$ with
$D_{\mathcal{B}}(Q)<\frac 12$. Fix five
elements $A_1,A_2,B_1,B_2,Y\in\mathcal{A}_{sa}$ that satisfy the conditions
(i)-(vii) from Lemma 5.1.

Notice that using (i), (iii), and (vii), it follows that
$q_{\mathcal{A}}(Y)=0$. Let us concentrate for the moment on the
element $Y\in\mathcal{A}_0$, which has $q_{\mathcal{A}_0}(Y)=0$.
On the one hand, $Y$ belongs to the thick AW*-subfactor $\mathcal{B}$.
On the other hand, by  condition (v), we have
$\mathbf{s}(Y)\leq Q$, so in particular we get
$D_{\mathcal{B}}\big(\mathbf{s}(Y)\big)<\frac 12$.
Using Theorem 5.1, there exists a spectrally symmetric element
$S\in\mathcal{B}_{sa}$, and a $2$-folding $\Phi=(Y_1,Y_2;S_1,S_2)$ in
$\mathcal{A}_0$, with
\begin{equation}
Y=Y_1+Y_2\text{ and }S=S_1+S_2.
\label{YS}\end{equation}

\begin{claim}
The elements $A_1,A_2,B_1,B_2,Y_1,Y_2,S_1,S_2$ all commute.
Moreover, one has the orthogonality relations
\begin{itemize}
\item[{\sc (a)}] $A_1,A_2,B_1\perp Y_1,Y_1,S_1,S_2$;
\item[{\sc (b)}] $B_2\perp S_1,S_2$.
\end{itemize}
\end{claim}
The relations {\sc (a)} are clear, since $A_1,A_2,B_1\perp
P$, and $Y_1,Y_2,S_1,S_2\in P\mathcal{A}P$.
To prove the relations {\sc (b)} we use the fact that
$S_1,S_2\in P\mathcal{A}P$, so that using condition (vi) from Lemma 5.1
and the fact that $\Phi=(Y_1,Y_2;S_1,S_2)$ is a folding, we have
$$B_2S_k=B_2PS_k=YS_k=(Y_1+Y_2)S_k=0,\,\,\,k=1,2.$$
Bsed on these orthogonality relations, we see that the only commutation that
needs to be checked is among $B_2,Y_1,Y_2$. Again using the fact that
$Y_1,Y_2\in P\mathcal{A}P$, and condition (vi) from Lemma 5.1, we have
$$B_2Y_k=B_2PY_k=YY_k=Y_kY=Y_kPB_2=Y_kB_2,\,\,\,k=1,2,$$
and we are done.

Having proven the Claim, we now define the elements
$X_1=(A_1-B_1)+(Y_1-S_1)$, $X_2=(A_2-B_2)+S$, and $X_3=Y_2-S_2$.
By the Claim, these three elements commute.
Moreover, $X_3$ is obviously spectrally symmetric. Using the orthogonality relations
from the Claim, combined with Lemma 5.1 and the features of
the $2$-folding $\Phi$, we also have
\begin{itemize}
\item $A_1-B_1$ and $Y_1-S_1$ are spectrally symmetric, and
$(A_1-B_1)\perp (Y_1-S_1)$,
\item $A_2-B_2$ and $S$ are spectrally symmetric, and
$(A_2-B_2)\perp S$,
\end{itemize}
hence $X_1$ and $X_2$ are also spectrally symmetric.

Finally, using condition (vii) from Lemma 5.1 and \eqref{YS}
we have
$$X_1+X_2+X_3=A_1-B_1+A_2-B_2+Y_1+Y_2-S_1-S_2+S=
A_1-B_1+A_2-B_2+Y=X,$$
and we are done.
\end{proof}

\begin{corollary}
Let $\mathcal{A}$ be an AW*-factor of type \twone, and let
$X\in\mathcal{A}_{sa}$ be an element with $q_{\mathcal{A}}(X)=0$. There exist
three commuting spectrally symmetric elements
$X_1,X_2,X_3\in \text{\rm Mat}_2(\mathcal{A})$ -- the $2\times 2$ matrix algebra --
such that
\begin{equation}
X_1+X_2+X_3=\left[\begin{array}{cc} X &0\\0&0\end{array}\right].
\label{cor2by2}\end{equation}
\end{corollary}
\noindent(According to Berberian's Theorem (see \cite{Be}), the matrix algebra
$\text{\rm Mat}_2(\mathcal{A})$ is an AW*-factor of type \twone.)
\begin{proof}
Denote the matrix algebra $\text{\rm Mat}_2(\mathcal{A})$ by
$\mathcal{A}_2$, and let $\tilde{X}\in \mathcal{A}_2$ denote the matrix in
the right hand side
of \eqref{cor2by2}. It is obvious that, if we consider the projection
$$
E=\left[\begin{array}{cc} 1 &0\\
0&0\end{array}\right],
$$
then $\mathbf{s}(\tilde{X})\leq E$. Since $D_{\mathcal{A}_2}(E)=\frac 12<1$,
and $q_{\mathcal{A}_2}(\tilde{X})=\frac 12 q_{\mathcal{A}}(X)=0$, the desired
conclusion follows imediately from Theorem 5.2.
\end{proof}

\end{document}